\documentclass[10pt]{article}   	% use "amsart" instead of "article" for AMSLaTeX format
\usepackage{geometry}                		% See geometry.pdf to learn the layout options. There are lots.
\usepackage{graphicx}				% Use pdf, png, jpg, or epsÂ§ with pdflatex; use eps in DVI mode
\usepackage{caption}
\usepackage{comment}
\usepackage{subcaption}
\usepackage{amssymb}
\usepackage{amsthm}
\usepackage{csquotes}
\usepackage{amsmath}
\usepackage{enumitem}
\usepackage{listings}
\usepackage{mathtools}
\usepackage{xcolor}
\usepackage{float}
\usepackage{physics}
\renewcommand{\thefootnote}{\arabic{footnote}}  % use numbers for footnotes instead of symbols
\makeatletter
\let\@fnsymbol\@arabic
\makeatother
\usepackage{aliascnt}
\usepackage[bookmarks=true,pdfstartview=FitH, pdfborder={0 0 0}, colorlinks=true,citecolor=red, linkcolor=blue]{hyperref}
\MakeOuterQuote{"}
\title{Stability of the Inviscid Power-Law Vortex}
\author{Tim Binz\thanks{Princeton University,
		Program in Applied \& Computational Mathematics, Fine Hall, Washington Road, 08544 Princeton, NJ, USA. Email: \url{tim.binz@princeton.edu}} , Matei P. Coiculescu\thanks{Princeton University,
		Department of Mathematics, Fine Hall, Washington Road, 08544 Princeton, NJ, USA. Email: \url{coiculescu@princeton.edu}}}

\theoremstyle{plain}
\newtheorem{thm}{Theorem}[section]
\newtheorem*{thm*}{Theorem}
\newtheorem*{conjecture}{Conjecture}
\newaliascnt{cor}{thm}
\newaliascnt{prop}{thm}
\newaliascnt{lem}{thm}
\newtheorem{cor}[cor]{Corollary}
\newtheorem{prop}[prop]{Proposition}
\newtheorem{lem}[lem]{Lemma}
\aliascntresetthe{cor}
\aliascntresetthe{prop}
\aliascntresetthe{lem} 
\theoremstyle{definition}
\newaliascnt{defn}{thm}
\newaliascnt{asu}{thm}
\newaliascnt{con}{thm}

\newtheorem*{defn*}{Definition}

\aliascntresetthe{defn}
\aliascntresetthe{asu}
\aliascntresetthe{con}
\newcounter{stp}
\newcounter{stpi}
\newcounter{stpci}
\newcounter{stpiii}

\theoremstyle{definition}
\newaliascnt{rem}{thm}
\newaliascnt{exa}{thm}
\newaliascnt{masu}{thm}
\newaliascnt{nota}{thm}
\newaliascnt{sett}{thm}
\aliascntresetthe{rem}
\aliascntresetthe{exa}
\aliascntresetthe{masu}
\aliascntresetthe{nota}
\aliascntresetthe{sett}
\newcommand\blfootnote[1]{
	\begingroup
	\renewcommand\thefootnote{}\footnote{#1}
	\addtocounter{footnote}{-1}
	\endgroup
}
\numberwithin{equation}{section}

\newcommand{\rL}{\mathrm{L}}

\newcommand{\R}{\mathbb{R}}
\newcommand{\N}{\mathbb{N}}

\newcommand{\Z}{\mathbb{Z}}

\begin{document}
	\maketitle

	\begin{abstract}		
		We prove that the power-law vortex $\overline{\omega}(x) = \beta \abs{x}^{-\alpha}$, which explicitly solves the stationary unforced incompressible Euler equations in $\mathbb{R}^2$ in both physical and self-similar coordinates, is exponentially linearly stable in self-similar coordinates with the natural scaling. This result, which is valid for functions in a weighted $\rL^2$ space and in the un-weighted $\rL^2$ space with a mild symmetry condition, answers a question from \cite{A}.
		Moreover, we prove that in physical coordinates the linearization around the power law vortex cannot generate an unstable $C_0$-semigroup.
	\end{abstract}
	
	\section{Introduction}
	\label{sec:intro}
	
	\blfootnote{2020 \emph{Mathematics Subject Classification:} 35Q31, 35B35, 47D06, 35C06
	}
	\blfootnote{Key Words and Phrases: Euler Equations, Linear Stability, Self-Similar Solutions}

		Consider the Cauchy problem of the two-dimensional Euler equations in vorticity form:
\begin{equation}
	\label{EULER}
	\left\{\begin{aligned}
		\partial_t \omega+ (v\cdot \nabla)\omega &=f\\
		\mathrm{K}_{\mathrm{BS}}\ast \omega(\cdot, t) &= v(\cdot, t) \\
		\omega(\cdot, t=0) &= \omega_{\mathrm{in}}
	\end{aligned}
	\right. 
\end{equation}
Here and in the sequel we denote the standard two-dimensional Biot-Savart kernel by $\mathrm{K}_{\mathrm{BS}}$, $f$ is a real-valued forcing term and the vorticity $\omega$ is a real-valued function defined on $\mathbb{R}^2\times [0,T)$. We consider solutions in the sense of distribution. In particular, we say $\omega$ is a (distributional) solution of the Euler  equation \eqref{EULER} if the following integral identity
\begin{equation}\label{DISTRO}
	\int_0^T \int_{\mathbb{R}^2} \big(\omega(\partial_t \phi + (\mathrm{K}_{\mathrm{BS}}\ast \omega)\cdot \nabla \phi) + f\phi \big)dxdt = -\int_{\mathbb{R}^2} \phi(x,0)\omega_{\mathrm{in}}(x)dx
\end{equation}
holds for every $\phi \in C^\infty_c(\mathbb{R}^2\times [0,T))$.
%Here and in the sequel we denote the standard two-dimensional Biot-Savart kernel by $\mathrm{K}_{\mathrm{BS}}$.
	
	\subsection{Well-posedness of the Euler Equations}
	
		We shall restrict ourselves to solutions of Equation \eqref{EULER} in a particular class of integrability. First, for any $1< q \leq \infty$ and $2<p\leq \infty$, we define:
	\begin{defn*}
		The function $\omega$ is in the class $\Upsilon_{q,p}^T$ if and only if
		$$\omega \in \rL^q_t([0,T], \rL^1_x\cap \rL^p_x) $$
		$$\mathrm{K}_{\mathrm{BS}}\ast \omega \in \rL^q_t([0,T], \rL^2_x) $$
	\end{defn*}
	\begin{defn*}
		The function $\omega$ is in the class $\Upsilon_p^0$ if and only if
		$$\omega \in \rL^1_x\cap \rL^p_x $$
		$$\mathrm{K}_{\mathrm{BS}}\ast \omega \in \rL^2_x $$
	\end{defn*}
	\begin{defn*}
		The function $\omega$ is in the class $\Upsilon_{q,p}^{\infty}$ if and only if $\omega\in \Upsilon_{q,p}^T$ for all $T\geq 0$.
	\end{defn*}
	We also use a slightly different notation if the domain of times does not include the time $t=0$:
	\begin{defn*}
		Let $a,b\in \mathbb{R}$. The function $\omega$ is in the class $\Upsilon_{q,p}^{[a,b]}$ if and only if
		$$\omega \in \rL^q_t([a,b], \rL^1_x\cap \rL^p_x) $$
		$$\mathrm{K}_{\mathrm{BS}}\ast \omega \in \rL^q_t([a,b], \rL^2_x).$$
		We also say that $\omega\in \Upsilon_{q,p}^{[a,\infty)}$ if and only if $\omega\in \Upsilon_{q,p}^{[a,b]}$ for every $b>a$.
	\end{defn*}
	One generally searches for solutions in $\Upsilon_{\infty,p}^\infty$ with initial data in $\Upsilon_p^0$. The famous theorem of Yudovich, proven in \cite{Y}, is:
	\begin{thm*}
		\label{thm:YUDO}
		Let $\omega_{\mathrm{in}}\in \Upsilon_\infty^0$ and let $f$ be some forcing term such that $f\in \Upsilon_{1,\infty}^\infty$. Then there exists a unique solution $\omega$ to Equation \eqref{EULER} in the class $\Upsilon_{\infty,\infty}^\infty$ with initial data $\omega_{\mathrm{in}}$.
	\end{thm*}
	
	It is a long-standing open question whether the Yudovich theorem can be extended to a broader class of solutions, in particular whether the Yudovich theorem holds in $\Upsilon_{\infty,p}^\infty$ for $p < \infty$.

	In a remarkable couplet of papers \cite{V1} and \cite{V2}, Vishik provides the first evidence that the so-called "Yudovich Class" $\Upsilon_\infty^0$ is sharp, proving:
	\begin{thm*}
		\label{thm:VISHIK}
		For every $2<p<\infty$, there exists $\omega_{\mathrm{in}}\in \Upsilon_p^0$ and a force $f\in \Upsilon_{1,p}^\infty$ with the property that there are uncountably many solutions $\omega\in \Upsilon_{\infty,p}^\infty$ to Equation \eqref{EULER} with initial data $\omega_{\mathrm{in}}$.
	\end{thm*}
	The monograph \cite{A} provides an alternative proof of Vishik's theorem while following a similar approach, and we shall generally use the notation and terminology from \cite{A}. Succinctly, one may describe Vishik's general strategy as constructing an unstable radial vortex in self-similar coordinates that generates non-uniqueness while breaking radial symmetry in physical coordinates. This follows the general program initiated by Jia and Sverak in \cite{JS1} and \cite{JS2} for proving non-uniqueness of Leray-Hopf solutions of the Naiver-Stokes equations. In \cite{JS1}, Jia and Sverak prove the striking result that under a spectral condition on the linearization of the Navier-Stokes equations in self-similar coordinates, one obtains at least two non-equal Leray-Hopf solutions. 
	
	\subsection{Self-similar Coordinates and Functional-analytic Set-up}
	\label{ssec:self-similar}
	
	We now discuss the interpretation of Vishik's proof in terms of dynamical systems proposed by the authors of the monograph \cite{A}. 	
	We first observe that Equation \eqref{EULER} admits many stationary solutions in the form of "radial vortices" or vortex profiles of the form:
	\begin{equation*} 
	\overline{\omega}(x) = g(\abs{x}), \quad \overline{v}(x) = \zeta(\abs{x}) x^\perp,
	\end{equation*} 
	where $x^\perp= (-x_2,x_1)$ and $\mathrm{K}_{\mathrm{BS}}\ast \overline{\omega} = \overline{v}$. 
	Suppose one could find a vortex profile $\overline{\omega}$ that is linearly unstable, for instance that one finds a real strictly positive eigenvalue $\lambda$ of the linearized Euler equations 
%	in self-similar coordinates 
	and a trajectory on the unstable manifold associated to $\lambda$ and $\overline{\omega}$ of the form
	\begin{equation}
	\omega = \overline{\omega} + \omega_{\mathrm{lin}} + o(e^{\lambda t}). \label{eq:omega construction}
	\end{equation} 
	Here $\omega_{\mathrm{lin}}= e^{\lambda t} \eta$ is a solution of the linearized Euler equations. We would then expect "non-uniqueness at time $t=-\infty$" because of the instability of the vortex. 
%	The outline proposed in \cite{A} is to instead 
	Recall that the solutions of the Euler equation \eqref{EULER} are invariant under the scaling
	\begin{equation*}
		\omega_{\mu,\lambda} (t,x)
		= \mu \cdot \omega(\mu t, \lambda x),
		\quad 
		v_{\mu,\lambda} (t,x)
		= \frac{\mu}{\lambda} \cdot 
		v (\mu t,\lambda x)
	\end{equation*}
	for $\mu, \lambda > 0$.
	The coordinates naturally adapted to this scaling are known as the forward self-similar coordinates of the Euler equation given by
	\begin{equation*}
		\begin{aligned}
			\xi &= xt^{-1/\alpha},  \quad &\tau &= \log(t), \\
			v(x,t) &= t^{1/\alpha-1}V(\xi,\tau), \quad &\omega(x,t) &= t^{-1}\Omega(\xi,\tau) 
		\end{aligned}
	\end{equation*}
	for a positive parameter $\alpha > 0$.
%	
%	$$\xi = xt^{-1/\alpha},  \quad \tau = \log(t)$$
%	$$v(x,t) = t^{1/\alpha-1}V(\xi,\tau), \quad \omega(x,t) = t^{-1}\Omega(\xi,\tau).$$
	Using these coordinates
	the Euler equations, without force, become
	\begin{equation*} 
		\begin{aligned}
	\partial_\tau \Omega -(1+\tfrac{\xi}{\alpha}\cdot\nabla_\xi)\Omega + V\cdot \nabla_\xi \Omega =0, \\
	V= \mathrm{K}_{\mathrm{BS}}\ast \Omega.
		\end{aligned} 
	\end{equation*} 
	If a self-similar profile $\Omega$ satisfies $\|\Omega\|_{\rL^p} = O(1)$ as $\tau \to -\infty$, then we also have that $\|\omega\|_{\rL^p}= O(t^{-1+\frac{2}{p\alpha}})$ as $t\to 0^+$. In particular, choosing $p=2/\alpha$ ensures that the Lebesgue norms are $O(1)$ in both coordinate systems, which is one reason why we consider this a "natural" choice of self-similar coordinates. To prove a non-uniqueness result in the spirit of Vishik, one should take $0<\alpha\leq 2/p$, which ensures the desired integrability. If one can find an unstable stationary solution $\overline{\Omega}$ of the self-similar equations, then one can hope to prove non-uniqueness at time $\tau=-\infty$, which corresponds to non-uniqueness at physical time $t=0$.
	Linearizing the Euler equation in self-similar coordinates around $\bar{\Omega}$ yields
	\begin{equation}
		\partial_\tau \Omega - L_{ss} \Omega = 0 ,
		\label{eq:linearization L_ss}
	\end{equation}
	where the linearization is given by
	\begin{equation}
		L_{ss}
		= (1+\tfrac{\xi}{\alpha}\cdot\nabla_\xi)\Omega - V\cdot \nabla_\xi \bar{\Omega} - \bar{V}\cdot \nabla_\xi \Omega .
		\label{eq:L_ss}
	\end{equation}
	Denoting an unstable eigenvalue by $\lambda$, the solution would appear like 
	\begin{equation}
		\Omega 
		= \bar{\Omega} + \Omega_{\mathrm{lin}} + \Omega_{\mathrm{per}}
		\label{eq:Omega construction}
	\end{equation}
	with $\Omega_{\mathrm{lin}} = e^{\lambda \tau} \eta$ solving the linearized Equation \eqref{eq:linearization L_ss}. Furthermore, 
	$\| \Omega_{\mathrm{per}} \| = o(e^{\lambda \tau})$ is obtained by a fixed point argument. Here $\bar{\Omega}$ corresponds to $\bar{\omega}$, $\Omega_{\mathrm{lin}}$ to $\|\omega_{\mathrm{lin}}\|= t^{\lambda}$, and $\Omega_{\mathrm{per}}$ yields the $o(t^\lambda)$-term. 
	In fact, without bothering with a perturbative corrective term to get a solution to the unforced equations, non-uniqueness of the Euler equations with a force follows as an immediate consequence of instability, see \cite{DM}. 
	
	\medskip 
	
	Let us know look closer into the functional-analytic setup of this strategy and how it is related to the integrability of the constructed (non-unique) solutions.
	While regularity and integrability of $\bar{\omega}$ is determined by the particular choice of the stationary solution, the regularity and integrability of $\omega_{\mathrm{lin}} + \omega_{\mathrm{per}}$, with $\| \omega_{\mathrm{per}} \| = o(t^\lambda)$, stems from the choice of the underlying Banach space $X$ that the linearization of the Euler equation in self-similar coordinates is realized in; the linearized operator is $L_{ss} \colon D(L_{ss}) \subset X \to X$ with domain $D(L_{ss}) = \{ \Omega \in X \colon L_{ss} \Omega \in X \}$.  
	Typical choices are weighted $\rL^q$-spaces given as $X = \rL^q(\R^2,\mu)$ for a measure $\mu(x) = |x|^\gamma \mathcal{L}^2(x)$ ($\mathcal{L}^2$ is the usual two-dimensional Lebesgue measure). 
	Since $\Omega_{\mathrm{lin}}$ evolves from an eigenfunction $\eta$ of $L_{ss}$ with an unstable eigenvalue $\lambda$, it lies in $D(L_{ss}^\infty)$ and may therefore have better regularity and integrability then provided in $X$. The bottleneck is the perturbative term $\Omega_{\mathrm{per}}$ that usually stems from a fixed point argument and therefore is usually in $X=\rL^q(\R^2,\mu)$.
	Using the embedding of the weighted $L^q$-spaces into weak Lebesgue spaces $X=\rL^q(\R^2,\mu) \hookrightarrow \rL^{r,\infty}(\R^2)$ for $r = \frac{2q}{2+\gamma}$ if $\gamma \geq 0$ we obtain that $\omega_{\mathrm{lin}} + \omega_{\mathrm{per}} \in \rL^{r,\infty}(\R^2)$. We now argue that any weighted $\rL^2$ space $X$ is inadequate for the completion of the program of Vishik if the background vortex is a radial power-law vortex. Indeed, suppose that $\overline{\omega}=|x|^{-\alpha} \in L^{2/\alpha,\infty}$. Since we expect that the perturbation is more regular than the background vorticity, we would need $\tfrac{2}{\alpha}<\frac{4}{2+\gamma}$ or $\alpha>1+\tfrac{\gamma}{2}\geq1$. However if $\alpha>1$ then the solution is not even locally in $\rL^{2}$.
	
	\smallskip 
	
	Note that the Biot-Savart operator $\omega\mapsto \mathrm{K}_{\mathrm{BS}}\ast \omega$ defined on Schwartz functions cannot be continuously extended to $\rL^2$, as shown in \cite{A}, it is well behaved on some closed linear subspaces of $\rL^2$. More generally, we shall show that for every $1<q <\infty$ there exist some closed linear subspaces of $\rL^q$, which we denote $\rL^q_m$, to which the operator can be continuously extended. Here $m\geq 2$ is an integer and $\rL^q_m$ is the space of $m$-fold rotationally symmetric functions lying in $\rL^q$. In other words, if $R_\theta: \mathbb{R}^2\to \mathbb{R}^2$ is the counterclockwise rotation of angle $\theta$ around the origin, then a function $f\in \rL^q_m$ satisfies
	$$f = f\circ R_{2\pi/m}.$$
	We observe that $L^2_1 = L^2$. Note, however, that we may still examine the spectral properties of the linearization of the Euler equations on spaces where $K_2\ast$ is unbounded. 
	
	\smallskip 
	
	Note that $L_{ss} \colon D(L_{ss}) \subset X \to X$ is a closed operator.
	We define its resolvent set $\rho(L_{ss},X)$ to be the open set of $z\in \mathbb{C}$ for which $L_{ss}-z$ has a bounded inverse from $X \to X$. We define the spectrum of $L_{ss}$, which we denote by $\sigma(L_{ss},X)$, to be the closed set which is the complement of the resolvent set. Since $L_{ss}$ depends on $\alpha$ we denote the spectrum by $\sigma_\alpha(L_{ss},X)$ if relevant. We shall often omit writing the dependence of the spectrum and resolvent set on the underlying Banach space whenever it is clear from the context. Instead, we may write $\sigma(L_{ss})$ or $\rho(L_{ss})$. The {\emph{spectral bound}} of $L_{ss}$ is defined to be
	$$s(L_{ss},X) := \inf_{\lambda\in \sigma(L_{ss},X)} \textrm{Re}(\lambda).$$
	From time to time we use $s_\alpha(L_{ss},X)$ to indicate the dependency on $\alpha$. In addition, we may simply denote the spectral bound by $s(L_{ss})$ if the underlying space is clear from the context.
	Recall that $\Omega_{\mathrm{lin}}$ solves
	\begin{equation*}
		\partial_{\tau} \Omega_{\mathrm{lin}} 
		- L_{ss} \Omega_{\mathrm{lin}} = 0,
		\quad \Omega_{\mathrm{lin}}(0) = \eta  
	\end{equation*}
	on $X$. This problem is well-posed if and only if $L_{ss}$ generates a strongly continuous semigroup $(T(t))_{t \geq 0}$ on $X$. The strongly continuous semigroup $(T(t))_{t \geq 0}$  satisfies
	\begin{equation*}
		\| T(t) \| \leq M_{\omega} e^{\omega t} 
	\end{equation*}
	for $M_{\omega} \geq 1$ and $\omega \in \R$. We call the semigroup \emph{quasi-contractive} if $M = 1$ and call
	\begin{equation*}
	\omega_0(L_{ss}) := \inf\{\omega\in \mathbb{R} : \exists M_\omega\geq 1 \textrm{ s.t. }\|\mathfrak{T}(t)\| \leq M_\omega e^{\omega t}, \quad \forall t\geq 0\}
	\end{equation*} 
	the \emph{growth bound} of $L_{ss}$.
%	Therefore the class of non-uniqueness constructed here is
%	\begin{equation*}
%		\rL^p(\R^2) + \rL^{m,\infty}(\R^2) .
%	\end{equation*}
%	Of course the most interesting case is when $p = m$. 

	\subsection{Main results}
	
	It is not difficult to see that the only stationary radial vortices solving the {\textbf{unforced}} self-similar Euler equations are precisely the power-law vortices of the form 
	$$\overline{\Omega}(\xi)= \beta(2-\alpha)\abs{\xi}^{-\alpha},$$
	and we note that these correspond exactly to radial power-law vortices solving the unforced stationary Euler equations in the original coordinates. The exponent $\alpha$ of the stationary profile is exactly determined by the choice of scaling for the self-similar coordinates, and the prefactor $\beta$ is an arbitrary real number. The natural question arises: are the power-law vortices unstable in the self-similar coordinates? This question was posed by the authors of \cite{A}. An affirmative answer would suggest that non-uniqueness for the unforced equations can arise from the (simple and explicit) power-law vortex, while a negative answer shows that a more complex self-similar profile that necessarily depends on the angular variable would have to be found if there is any hope to complete the program proposed in \cite{A}. 
	Loosely speaking our main results show that the answer to this question is negative. Therefore, our results rule out the power-law vortex as a background solution in the approach of Vishik. This confirms the statement in~\cite{ABC} that finding such an unstable vortex is "far from elementary".
%	In order to give a more precise statement we have to introduce the Banach spaces we are working in. 
	Since the power-law vortex is singular, it is natural to work in a weighted $\rL^2$ space:
	\begin{equation*} 
		\rL^2_w(\mathbb{R}^2) := \rL^2(\mathbb{R}^2, |x|^{2+\alpha}dx).
	\end{equation*} 
	Following Vishik, we consider the subspace of $m$-fold symmetric vector-fields $\rL^2_{w,m}(\R^2)$.
	The weighted $\rL^2$-based space $\rL^2_{w,m}$ is in fact the unique such space for which the nonlocal part of the linearization $\Omega\to (K_2\ast \Omega) \cdot\nabla \overline{\Omega}$ becomes a symmetric operator. This parallels similar choices made when studying the Euler equations in self-similar variables, as for example, in \cite{TB}.

	We are now prepared to state our first main result:
	
	\begin{thm}
		\label{thm:MAINMAIN}
		Let $\beta\in \mathbb{R}$. Let $m\geq 2$ and $\alpha \in (0,2)$ or $m=1$ and $\alpha \in (0,2-\sqrt{2})$.  
		If $\overline{\Omega}(\xi) = \beta(2-\alpha) \abs{\xi}^{-\alpha}$ is the radial power-law vortex that solves the unforced Euler equations in self-similar coordinates and $L_{ss}$ is the linearization of the Euler equations around $\overline{\Omega}$, then the operator $L_{ss}$ generates a quasi-contractive semigroup on $\rL^2_{w,m}(\R^2)$ with growth bound
		\begin{equation*}
			\omega_0(L_{ss}) = s_\alpha(L_{ss}) = \frac{1}{2}-\frac{2}{\alpha} < 0.
		\end{equation*}
	\end{thm}
	
	Note that our upper bound for the spectral bound is always independent of the "size" of the background vortex, which is controlled by the parameter $\beta$, so the size of the background vortex cannot contribute to an instability (in contrast with the program of Jia and Sverak).

	The weighted space excludes perturbations that are in $\rL^2$ or that are even more singular near the spatial origin. The weighted space, however, does not exclude a perturbation in $\rL^2$ that does not decay sufficiently quickly at infinity. Considering that Vishik's example arises from a local singularity at the spatial origin and that Vishik's eigenfunction lies in our weighted space, we consider this as an inessential deficiency.
	
	\smallskip 
	
	Let us now discuss the role of the parameter $\beta$. In order to distinguish between power-law vortices for different $\beta$ we write $\bar{\Omega}_\beta$. 
	Note that the linearization can be written as
	\begin{align*}
		L_{ss,\beta}
		= (1+\tfrac{\xi}{\alpha}\cdot\nabla_\xi)\Omega - V\cdot \nabla_\xi \bar{\Omega}_\beta - \bar{V}_\beta \cdot \nabla_\xi \Omega 
		= \beta \cdot 
		\biggl(
		\beta^{-1} (1+\tfrac{\xi}{\alpha}\cdot\nabla_\xi)\Omega
		- V\cdot \nabla_\xi \bar{\Omega}_1 - \bar{V}_1 \cdot \nabla_\xi \Omega 
		\biggr) .
	\end{align*}
	The term in the parentheses consists of a term of order $\beta^{-1}$ and the linearization of the Euler equation in physical coordinates around $\bar{\Omega}_1$. We denote by $L_{\phi} \colon D(L_\phi) \subset X \to X$ the linearization given by
	\begin{equation*}
		L_\phi \Omega = 
		- V\cdot \nabla_\xi \bar{\Omega}- \bar{V} \cdot \nabla_\xi \Omega ,
		\quad D(L_\phi) = \{ \omega \in X \colon L_{\phi} \omega \in X \}
	\end{equation*}
	with $V = K_{\mathrm{BS}} \ast \Omega$.
	Formally, we have that the term in the parentheses satisfies
	\begin{equation*}
		\biggl(
		\beta^{-1} (1+\tfrac{\xi}{\alpha}\cdot\nabla_\xi)\Omega
		- V\cdot \nabla_\xi \bar{\Omega}_1 - \bar{V}_1 \cdot \nabla_\xi \Omega 
		\biggr)
		\to L_\phi \Omega \quad \text{ as } \beta \to + \infty . 
	\end{equation*} 
	Despite the formal limit, this operation may not be well-defined since $D(L_{ss}) \cap D(L_\phi) = \{ 0 \}$ can occur. However, we can make the convergence precise in the resolvent sense. This formal limit was used in \cite{ABC} to show the instability of $L_{ss}$ for $\beta > 0$ sufficiently large from an instability for $L_\phi$.
	Since our main result establishes the stability of 
	$L_{ss}$, it is natural to ask whether the same property holds for 
	$L_{\phi}$.
%	We also prove a result on the linearization of the Euler equations around the power-law vortex in physical coordinates, which we shall write as:
%	$$(\partial_t - L_\phi)\Omega=0. $$ 
	Due to the time-irreversibility of the Euler equations, the spectrum of $L_\phi$ is symmetric about the imaginary axis. We shall prove that the spectrum lies in the imaginary axis:
	\begin{thm}
		\label{thm:physicalMAINMAIN}
		Let $\beta\in \mathbb{R}\setminus\{0\}$. Let $m\geq 2$ and $\alpha \in (0,2)$ or $m=1$ and $\alpha \in (0,2-\sqrt{2})$. 
		Let $\overline{\Omega}(\xi) = \beta(2-\alpha) \abs{\xi}^{-\alpha}$ be the radial power-law vortex that solves the unforced Euler equations in self-similar coordinates and let $L_{\phi}$ be the linearization of the Euler equations around $\overline{\Omega}$. Then the operator $L_{\phi}$ generates a continuous group of isometries on $\rL^2_{w,m}(\R^2)$ with
		$$\sigma(L_\phi) \subset \{z \in \mathbb{C} : \textrm{Re}(z)=0\}.$$
	\end{thm}
	
	Now, with \autoref{thm:physicalMAINMAIN} in hand, we may view our results in the context of classical results on hydrodynamics stability. In the comprehensive work \cite{C} by Chandrasekhar, various stability and instability results for fluid motion are described, including the well-known criterion of Rayleigh's for stability of steady inviscid flow between two co-axial cylinders. Rayleigh's criterion would suggest, but not rigorously prove in the case of our infinite energy power-law vortices on the whole plane, that a vorticity profile $\abs{x}^{-\alpha}$ is stable in physical coordinates whenever $0<\alpha<2$. We consider our \autoref{thm:physicalMAINMAIN} as the first mathematically rigorous proof of this expected stability. We also mention the more recent work \cite{ZZ} of Zelati and Zillinger, in which the authors consider the linear stability of vorticity profiles with singularities of power-law type in physical coordinates. The work in \cite{ZZ}, however, does not provide quantitative information on the spectrum and does not consider the stability in self-similar coordinates, which is the main part of the present article.
	
	\smallskip 
	
	We now remark the relationship our work has with the stability theorems proven by Arnold in \cite{AR} and the related questions proposed by Yudovich in \cite{Y2}. Let $\psi$ be the stream function of a stationary solution of the incompressible Euler equations. Arnold proves in \cite{AR} that the stationary flow is stable if the velocity profile is convex, or $\nabla \psi/ \nabla \Delta \psi>0$. When the velocity profile is concave, or $\nabla \psi/ \nabla \Delta \psi<0$, then there are finitely many unstable eigenvalues of the corresponding linear problem. The velocity profile corresponding to the power-law vortex $\overline{\Omega}$ is concave, so we expect finitely many unstable eigenvalues. Our work improves this claim to a rigorously proven statement that there are no unstable eigenvalues for the linearization around $\overline{\Omega}$. Yudovich in \cite{Y2} proposes that understanding the stability (or instability) of ideal fluid flows is an important problem in mathematical hydrodynamics, and we consider our work as progress in that direction.
	
	\subsection{\texorpdfstring{Complementary Results in $\rL^p$-spaces}{Complementary results in rLp-spaces}}
	
	Despite being less natural than the weighted case, we proceed by making a few remarks about the unweighted setting. The following result is an analogue of our main \autoref{thm:MAINMAIN}:
	\begin{prop}
		\label{BUSYPROP}
		Let $\beta\in \mathbb{R}$ and $\alpha \in (0,1)$.
		Let $m\geq 3$ be any integer. If $\overline{\Omega}(\xi) = \beta(2-\alpha) \abs{\xi}^{-\alpha}$ is the radial power-law vortex that solves the unforced Euler equations in self-similar coordinates and $L_{ss}$ is the linearization of the Euler equations around $\overline{\Omega}$, then the point spectrum of $L_{ss}$ is empty and
		$$1-\frac{1}{\alpha}\leq s_\alpha(L_{ss},\rL^2_m ) \leq 1-\frac{1}{\alpha } + \frac{4}{\alpha (m-2)}.$$
	\end{prop}
	Therefore, the operator $L_{ss}$ is spectrally stable with bound $\leq 0$ whenever
	$$m \geq \frac{6-2\alpha}{1-\alpha}.$$
	Note that in Vishik's work \cite{V1,V2} he also uses some (probably large) $m$. 
	Our method in attaining this spectral bound, which we elaborate on in the \autoref{sec:Surjection}, is certainly not optimal. We conjecture that our result can be sharpened as follows:
	\begin{conjecture}
		Let $\beta\in \mathbb{R}$. Let $m\geq 3$ and $\alpha \in (0,1)$.   If $\overline{\Omega}(\xi) = \beta(2-\alpha) \abs{\xi}^{-\alpha}$ is the radial power-law vortex that solves the unforced Euler equations in self-similar coordinates and $L_{ss}$ is the linearization of the Euler equations around $\overline{\Omega}$, then the operator $L_{ss}$ generates a quasi-contractive semigroup on $\rL^2_m(\R^2)$ with growth bound
		\begin{equation*}
			\omega_0(L_{ss}) = s_\alpha(L_{ss}) = 1-\frac{1}{\alpha} < 0.
		\end{equation*}
	\end{conjecture}
	As discussed in \autoref{ssec:self-similar}, the unweighted $\rL^2$-case might not be the most natural setting for the linear analysis. We therefore continue with remarks about the case when $X=\rL^q$ and $q\neq 2$. 
	The technical reason for working with $\rL^2$ is that
	by the Plancherel theorem we have
	$$\rL^2_m = \bigoplus_{k\in \mathbb{Z}}^{\ell^2} U_{km},$$
	where 
	$$U_{km} := \{f(r)e^{imk\theta} : f\in \rL^2(\mathbb{R}^+, rdr)\}$$
	are mutually orthogonal subspaces. 
	This decomposition can never occur when $q \neq 2$. Nevertheless we may consider the operator $L_{ss}$ restricted to the closed, proper subspaces $U_{k,q}$ of $\rL^q$ given by
	$$U_{k,q} := \{f \colon f(r,\theta) = f(r)e^{ik\theta} : f\in \rL^q(\mathbb{R}^+, rdr)\}.$$
	Here we are able to prove a weaker result: an upper bound on the spectral bound.
	%In the case when $q=2$, we use the additional Hilbert space structure and are able to prove a much stronger result on the whole of $\rL^2_m$. 
%	More precisely, for the Banach case we prove:
	
	\begin{prop}
		\label{thm:MAINMAINBANACH}
		Let $\beta\in \mathbb{R}$ and $\alpha \in (0,1)$.
		Let $q\neq 2$ be such that $1< q \leq 2/\alpha$ and let $k$ be any integer satisfying $\abs{k}\geq 3$. If $\overline{\Omega}(\xi) = \beta(2-\alpha) \abs{\xi}^{-\alpha}$ is the radial power-law vortex that solves the unforced Euler equations in self-similar coordinates and $L_{ss}$ is the linearization of the Euler equations around $\overline{\Omega}$, then 
		$$s_\alpha(L_{ss}, U_{k,q} ) \leq 1-\frac{2}{\alpha q} + \frac{8}{\alpha q(\abs{k}-2)}.$$
	\end{prop}
	
	%In the Hilbert space case, we shall prove the %following stronger theorem:
	%\begin{thm}
	%\label{thm:MAINMAIN}
	%Let $\beta\in \mathbb{R}$ and $\alpha \in (0,1)$. Let $m\geq 3$ be any integer. If $\overline{\Omega}%(\xi) = \beta(2-\alpha) \abs{\xi}^{-\alpha}$ is %the radial power-law vortex that solves the %unforced Euler equations in self-similar %coordinates and $L_{ss}$ is the linearization of %the Euler equations around $\overline{\Omega}$, %then the operator $L_{ss}$ generates a quasi-%contractive semigroup on $\rL^2_m(\R^2)$ with
	%    \begin{equation*}
		%        \omega_0(L_{ss}) = s(L_{ss}) = 1-\frac{1}%{\alpha} < 0 .
		%   \end{equation*}
	%\end{thm}
	
	Before continuing, we offer a few remarks. First, we have spectral linear stability (a spectral bound less than zero) on $U_{k,q}$ if 
	$$2+\frac{8}{2-\alpha q}<\abs{k}.$$
	Second, it is not obvious from the resolvent bound we find in \autoref{BUSYPROP} and \autoref{thm:MAINMAINBANACH} that $L_{ss}$ is an operator on $L^2_m$ and respectively $U_{k,q}$ satisfying the hypotheses of the Hille-Yosida  theorem. Thus, we cannot immediately conclude that $L_{ss}$ generates a continuous semigroup in the non-weighted $\rL^2$ case or the non-Hilbert space case (note also that the growth bound is not a priori equal to the spectral bound). %Third, our upper bound for the spectral bound is always independent of the "size" of the background vortex, which is controlled by the parameter $\beta$, so the size of the background vortex cannot contribute to an instability (in contrast with the program in Jia and Sverak).
	We finish this section by an analogous result to \autoref{thm:physicalMAINMAIN} about the linearization of the Euler equation in physical coordinates. 
	
	\begin{prop}
		\label{prop:physicalMAINMAINunweighted}
		Let $\beta\in \mathbb{R}\setminus\{0\}$. Let $m\geq 1$ and $\alpha \in (0,1)$. % or $m=1$ and $\alpha \in (0,2-\sqrt{2})$. 
		Let $\overline{\Omega}(\xi) = \beta(2-\alpha) \abs{\xi}^{-\alpha}$ be the radial power-law vortex that solves the unforced Euler equations in self-similar coordinates and let $L_{\phi}$ be the linearization of the Euler equations around $\overline{\Omega}$. If $L_{\phi}$ generates a $C_0$-semigroup on $\rL^2_{m}(\R^2)$ then
		$$\sigma(L_\phi) \subset \{z \in \mathbb{C} : \mathrm{Re}(z)=0\}.$$
	\end{prop}
	
	Let us point out that the proof of  \autoref{prop:physicalMAINMAINunweighted} does not really use the fact that we work in $\rL^2_m(\R^2)$. It essentially relies on the fact that the power law vortex is $\alpha$-homogeneous. Therefore it seems to be generalizable to other Banach spaces. 
	
	\bigskip 
	
		Besides the application to fluid dynamics, we consider an intriguing aspect of our work to be the novel techniques we use in our analysis of the Euler equations linearized around the singular power-law vortex. The singularity of the background vortex leads to unbounded operators $L_{ss}, L_\phi$ with unbounded coefficients, which, to our knowledge, cannot be handled easily by any previously known method. For example, $L_{ss}$ cannot be thought of as the compact perturbation of a skew-adjoint operator, which would be case if the background vortex were smooth. Our techniques might be of use to other linear problems with singular coefficients or to other nonlinear stability problems around singular background solutions.
	
	\subsection{Outline of the article}
	
	The article is organized as follows: 
	In \autoref{sec:power law vortex} we prove the fact that the power law vortices are the only stationary and self-similar solutions of the Euler equation. 
	In \autoref{sec:BSLAW} we show that the Biot-Savart operator can be extended continuously to $\rL^q_m(\R^2)$. 
	With this preparation we prove our main \autoref{thm:MAINMAIN} in \autoref{sec:linear stability self-similar}. Finally, in \autoref{sec:linear stability} we show \autoref{thm:physicalMAINMAIN}.
	Additionally, we prove the results \autoref{BUSYPROP} and \autoref{thm:MAINMAINBANACH} in \autoref{sec:Surjection}. In \autoref{sec:explicit solution} we show that the point spectrum of $L_{ss}$ is empty. In \autoref{sec:mathematica}, we provide some Mathematica code used in previous sections. Finally, in \autoref{sec:physicalunweighted} we show \autoref{prop:physicalMAINMAINunweighted}.
	
	\subsection*{Acknowledgments}
	We thank Princeton University for its support. 
	MPC thanks Professor Camillo De Lellis for his constant encouragement and support during MPC's graduate studies, as well as for sharing his deep insight into this and many other problems of mathematics. MPC acknowledges the support of the National Science Foundation in the form of an NSF Graduate Research Fellowship.
	TB acknowledges the support of the DFG Walter-Benjamin Fellowship no. 
	538212014.	

	\section{The Inviscid Power-Law Vortex}
	\label{sec:power law vortex}
	
		We now more precisely discuss the change of coordinates used by Vishik in \cite{V1} and \cite{V2} as well as by the authors of \cite{A}. Given a solution $\omega$ of Equation \eqref{EULER} on $\mathbb{R}^2 \times [T_0,T_1]$, we introduce a function $\Omega : \mathbb{R}^2\times [\log{T_0}, \log{T_1}]\to \mathbb{R}$ given by the following transformation. We set $\tau = \log{t}$, $\xi = xt^{-1/\alpha}$, and let 
	$$\Omega(\xi, \tau) = e^\tau \omega(e^{\tau/\alpha} \xi, e^\tau).$$
	The reverse transformation is given by:
	$$\omega(x,t) = t^{-1} \Omega(t^{-1/\alpha} x, \log{t}).$$
	If the vector field $v$ is given by $(\mathrm{K}_{\mathrm{BS}}\ast \omega)$ and the vector field $V$ is given by $(\mathrm{K}_{\mathrm{BS}}\ast \Omega)$, then we have the following transformation rules:
	$$V(\xi, \tau) = e^{\tau(1-1/\alpha)}v(e^{\tau/\alpha}\xi, e^\tau),$$
	$$v(x,t) = t^{-1+1/\alpha} V(t^{-1/\alpha} x, \log{t}).$$
	Vishik in \cite{V1} and \cite{V2} as well as the authors of \cite{A} show that if $\omega$ satisfies Equation \eqref{EULER}, then the function $\Omega$ satisfies 
	\begin{equation}
		\label{SEULER}\left\{
		\begin{aligned}
			\partial_\tau \Omega - \bigg(1+ \frac{\xi}{\alpha}\cdot \nabla\bigg)\Omega+ (V\cdot \nabla)\Omega &=0\\
			\mathrm{K}_{\mathrm{BS}}\ast \Omega(\cdot, \tau) &= V(\cdot, \tau). \\
		\end{aligned}\right.
	\end{equation}
	We usually refer to this system of coordinates as the exponential self-similar coordinates or simply the self-similar coordinates. We refer to the system \eqref{SEULER} as the vorticity form of the Euler equations in self-similar coordinates.

	Now we introduce exponential self-similar polar coordinates. In particular, if $(r,\theta, t)$ are the usual polar coordinates on $\mathbb{R}^2\times [T_1,T_2]$, we let $\rho = rt^{-1/\alpha}$, $\theta$ unchanged, and $\tau = \log{t}$ to get a new $(\rho,\theta,\tau)$ system of coordinates. We consider the velocity form of the Euler Equations:
	\begin{equation}
		\label{VEULER}
		\begin{gathered}
			\partial_t v+ (v\cdot \nabla)v =-\nabla p\\
			\textrm{div } v=0
		\end{gathered}
	\end{equation}
	which we write in terms of exponential self-similar polar coordinates as:
	\begin{equation}
		\label{VPOLAR}\left\{
		\begin{aligned}
			\partial_\tau V_\rho +\big(\frac{1}{\alpha}-1\big)V_\rho - \frac{\rho}{\alpha}\partial_\rho V_\rho +V_\rho \partial_\rho V_\rho +\frac{V_\theta}{\rho}\partial_\theta V_\rho -\frac{V_\theta^2}{\rho} = - \partial_\rho P\\
			\partial_\tau V_\theta +\big(\frac{1}{\alpha}-1\big)V_\theta - \frac{\rho}{\alpha}\partial_\rho V_\theta +V_\rho \partial_\rho V_\theta +\frac{V_\theta}{\rho}\partial_\theta V_\theta +\frac{V_\rho V_\theta}{\rho} = - \frac{1}{\rho}\partial_\theta P\\
			\partial_\rho(\rho V_\rho) + \partial_\theta V_\theta =0.
		\end{aligned}\right.
	\end{equation}
	
	One can see that, up to a constant prefactor, the only stationary radial vortex $V=V_\theta(\rho)e_\theta$ satisfying the exponential self-similar equations is precisely the power-law vortex with $V_\theta(\rho) = \beta\rho^{1-\alpha}$, where $\beta$ is any real number. 
	
	\begin{prop}
		The unique solution of Equation \eqref{VPOLAR} of the form $V= V_\theta(\rho)e_\theta$ is given by the profile $V_\theta(\rho) = \beta\rho^{1-\alpha}$, where $\beta$ is any real number. 
	\end{prop}
	\begin{proof}
		The divergence-free condition is clearly satisfied. The second listed equation in Equation \eqref{VPOLAR} simplifies to $(1/\alpha-1)V_\theta(\rho)-(p/\alpha)V_\theta'(\rho)=0,$ a first-order differential equation whose unique solution is $V_\theta(\rho)=\beta \rho^{1-\alpha}$. For an appropriate choice of the pressure $P$, the equation listed first in Equation \eqref{VPOLAR} will also be satisfied.
	\end{proof}
	
	We are thus led to define the velocity profile:
	\begin{equation}
		\overline{V}(\rho) =  \beta\rho^{1-\alpha} e_\theta
	\end{equation}
	and its associated vorticity profile:
	\begin{equation}
		\overline{\Omega}(\rho) =\beta(2-\alpha)\rho^{-\alpha}.
	\end{equation}
	
	Note that $\overline{\Omega}$ corresponds exactly by the transformation back to physical coordinates to the stationary radial vortex with vorticity profile $\beta(2-\alpha)r^{-\alpha}$.

\section{The Biot-Savart Law}
\label{sec:BSLAW}

We now recall some useful notation. We denote $\rL^q_m$ to be the closed linear subspace of $\rL^q$ of elements that are $m$-fold symmetric. In other words, if $R_\theta: \mathbb{R}^2\to \mathbb{R}^2$ is the counterclockwise rotation of angle $\theta$ around the origin, then a function $f\in \rL^q_m$ satisfies
	$$f = f\circ R_{2\pi/m}.$$
	It will be convenient to define the following closed linear subspaces of $\rL^q_m$:
	$$U_{km,q} := \{f(r)e^{imk\theta} : f\in \rL^q(\mathbb{R}^+, rdr)\}.$$
	In the case when $q=2$, the subspaces $U_{km,2}$ are mutually orthogonal, and we have the direct sum:
	$$\rL^2_m = \bigoplus_{k\in \mathbb{Z}}^{\ell^2} U_{km,2}.$$
	The topology of the space on the right is given by the $\ell^2$ direct sum of the spaces $U_{km,2}$ with the norm of an element in the direct sum $\sum_{k\in\mathbb{Z}} f_k$ given by
	$$\bigg(\sum_{k\in\mathbb{Z}} \|f_k\|_{\rL^2}^2\bigg)^{1/2}.$$
	By the Theorem of Plancherel, the direct sum with this topology is exactly the space $\rL^2_m$ with the subspace topology inherited from $\rL^2$. In the case when $q\neq 2$, we no longer have the benefit of the Hilbert space structure. Indeed, it is not difficult to see that $\rL^q_m$ cannot be written as the $\ell^p$ direct sum of the spaces $U_{km,q}$ for any $p\in [1,\infty]$, see \cite{E1} and Section 13.5 on the Hausdorff-Young inequality in \cite{E2}. Instead, we shall work in each closed proper subspace $U_{k,q}$ of $\rL^q_m$ separately. 
	We denote the Schwartz space by $\mathcal{S}$ and its dual by $\mathcal{S}^*$. We begin by stating a slight improvement of a lemma from \cite{A}. The improvement allows for a continuous extension of $\mathrm{K}_{\mathrm{BS}}\ast$ to $\rL^q_m$ whenever $1< q<\infty$.
	
	\begin{lem}
		\label{lem:K2}
		For every $m\geq 2$ there exists a unique continuous operator $T: \rL^q_m \to S^\ast$ satisfying:
		\begin{enumerate}
			\item If $\varphi\in \mathcal{S}$, then $T\varphi = \mathrm{K}_{\mathrm{BS}}\ast \varphi$ in the sense of distribution.
			\item There exists a constant $C>0$ such that for every $\varphi \in \rL^q_m$ there exists $v(\varphi):= v\in W^{1,q}_{loc}$ such that
			\begin{itemize}
				\item $R^{-1}\|v\|_{\rL^q(B_R)}+ \|Dv\|_{\rL^q(B_R)}\leq C \|\varphi\|_{\rL^q(\mathbb{R}^2)}$ for all $R>0$
				\item $\textrm{div } v=0$ and $T(\varphi)=v(\varphi)$ in the sense of distribution.
			\end{itemize}
		\end{enumerate}
	\end{lem}
	Certainly, since $T$ is a continuous linear operator on $\rL^q_m$, $T$ is a continuous linear operator when restricted to the closed invariant subspace $U_{k,q}\subset \rL^q_m$.
	We shall denote $T$ interchangeably with $\nabla^\perp\Delta^{-1}$ or $\mathrm{K}_{\mathrm{BS}}\ast$, although we shall mostly use the latter notation.

	\subsection{Proof of \autoref{lem:K2}}
    \begin{proof}
	Let $f\in \rL^q_m \cap \mathcal{S}$ and let $v= \mathrm{K}_{\mathrm{BS}} \ast f$. We first claim that 
	\begin{equation}
		\label{ZEROINTEGRAL}
		\int_{B_R} v =0\quad \textrm{ for every } R>0.\end{equation}
	Indeed, we have $v= \nabla^\perp h$, where $h$ is the unique classical solution of $\Delta h = f$ given by $K\ast f$ where $K(x) = \frac{1}{2\pi}\log\abs{x}$. Since the kernel $K$ is invariant under all rotations $R_\theta$ and since $f$ is $m$-fold symmetric, we conclude that $h$ is $m$-fold symmetric. Therefore, 
	$$R_{-2\pi/m}\nabla h(R_{2\pi/m} x)= \nabla h(x).$$
	Thus, integrating in $x$ we conclude 
	$$\int_{B_R} \nabla h = R_{2\pi/m} \int_{B_R} \nabla h,$$
	so 
	$$\int_{B_R} \nabla h = \frac{1}{m}\sum_{k=0}^{m-1} R_{2k\pi/m} \int_{B_R} \nabla h.$$
	However, since $m\geq 2$, the sum is zero, which shows that $\int_{B_R} \nabla h=0$. Finally, with the property just shown, we may use the Poincar\'e inequality to conclude:
	$$R^{-1}\|v\|_{\rL^q(B_R)} + \|Dv\|_{\rL^q(B_R)} \leq C \|f\|_{\rL^q(B_R)}$$
	since $\|Dv\|_{\rL^q}\leq \|f\|_{\rL^q}$ by the Calder\'on-Zygmund theorem.
\end{proof}

	\section{Linear Stability of the Power-Law Vortex in Self-Similar Coordinates}
	\label{sec:linear stability self-similar}
The goal of this section is to prove the following theorem, which is the final result of a highly technical analysis of the linearization of the Euler equations around a singular power-law vortex. 
\begin{thm}
		\label{thm:mainss}
		Let $\beta\in \mathbb{R}$. Let $m\geq 2$ and $\alpha \in (0,2)$ or $m=1$ and $\alpha \in (0,2-\sqrt{2})$.  
		If $\overline{\Omega}(\xi) = \beta(2-\alpha) \abs{\xi}^{-\alpha}$ is the radial power-law vortex that solves the unforced Euler equations in self-similar coordinates and $L_{ss}$ is the linearization of the Euler equations around $\overline{\Omega}$, then the operator $L_{ss}$ generates a quasi-contractive semigroup on $\rL^2_{w,m}(\R^2)$ with growth bound
		\begin{equation*}
			\omega_0(L_{ss}) = s_\alpha(L_{ss}) = \frac{1}{2}-\frac{2}{\alpha} < 0.
            	\end{equation*}
	\end{thm}
We write the linearization around $\overline{\Omega}$ in self-similar coordinates as
	$$(\partial_\tau -L_{ss})\Omega=0.$$
	The main technical issue, for which no standard method can be used, is the highly singular behavior of the coefficients of $L_{ss}$. For example, we cannot simply say that $L_{ss}$ is the compact perturbation of a skew-adjoint operator, which, on the contrary, would be the case if the background vortex profile $\overline{\Omega}$ were smooth.

	\subsection{Preliminaries}
	\label{ssec:prelims}

	As before, we attempt to reduce our problem to the analysis of an infinite dimensional collection of one-dimensional problems by using the invariance of the operator with respect to certain subspaces. 
	By Plancherel's theorem we have
	\begin{equation}
		\rL^2_{w,m} = \bigoplus^{\ell^2}_{k\in\mathbb{Z}} U_{km}  ,  
		\label{eq:plancherel}
	\end{equation}
	where $U_n$ denotes the proper subspace of $\rL^2_w(\R^2)$ given by
	\begin{equation*}
		U_n := \{ f \colon f(r,\theta) = g(r)e^{in\theta} \text{ and } g \in \rL^2(\R_+,r^{3+\alpha}\mathrm{d}r) \} . 
	\end{equation*}
	The equation \eqref{eq:plancherel} should be understood as equality of Hilbert spaces, i.e. the scalar product satisfies
	\begin{equation*}
		\langle U,V \rangle_{\rL^2_m}
		= \sum_{k \in \Z} \langle f_k, g_k \rangle_{\rL^2(\R,r^{3+\alpha} \mathrm{d}r)},
	\end{equation*}
	where $f_n,g_n$ denotes the component of $U,V$ in $U_{nm}$. Here
	$\langle \, \cdot \,, \, \cdot \, \rangle$ denotes the usual $\rL^2$-pairing.
    
	In what follows, we shall abuse notation and denote the (exponential self-similar) radial coordinate by $r$ instead of $\rho$.
	We recall that the inverse Laplacian $\Delta^{-1}$ preserves rotational symmetries. Thus, for any $m$-fold symmetric $\Omega$, $\Delta^{-1}\Omega$ is also $m$-fold symmetric. We formally decompose:
	\begin{equation}\label{DECOMPOSE}\Delta^{-1} \Omega = \sum_{k\in \mathbb{Z}} f_k(r)e^{mik\theta}\end{equation}
	where $f_k(r) \in \rL^2(\mathbb{R}^+, r^{3+\alpha}\mathrm{d}r)$ for all $k\in \mathbb{Z}$. %Let $\lambda= \lambda_1+i\lambda_2$ be our putative element of the spectrum, where $\lambda_1,\lambda_2\in \mathbb{R}$ and $\lambda_1>1-\frac{2}{q\alpha}$. We consider the eigenvalue equation $L_{ss}\Omega-\lambda\Omega=0$, or equivalently
	Using the decomposition \eqref{DECOMPOSE} the operator
    \begin{equation*}
        L_{ss}\Omega 
        = \Omega +\big(\tfrac{\xi}{\alpha}-\overline{V}\big)\cdot\nabla\Omega-\nabla^{\perp}\Delta^{-1}\Omega\cdot\nabla\overline{\Omega}
    \end{equation*}
	%, the above linear partial differential equation 
    becomes equivalent to an infinite family of ordinary differential operators. We have 
    %The ordinary differential operators corresponding to the parameter $k$ are:
	\begin{align*}
    (L_{ss} \Omega)(r,\theta) 
    &= \sum_{k \in \Z} \biggl( 
    (\partial_r^2+\tfrac{\partial_r}{r}+\tfrac{\partial_\theta^2}{r^2})(f_k(r)e^{mik\theta}) \\
    &+\big(\tfrac{r}{\alpha}e_r-\beta r^{1-\alpha}e_\theta\big)\cdot(\partial_r\big((\partial_r^2+\tfrac{\partial_r}{r}+\tfrac{\partial_\theta^2}{r^2})(f_k(r)e^{mik\theta})\big)e_r \\
    &+\tfrac{\partial_\theta}{r}\big((\partial_r^2+\tfrac{\partial_r}{r}+\tfrac{\partial_\theta^2}{r^2})(f_k(r)e^{mik\theta})e_\theta\big)
    \\
    &+ \big(\tfrac{\partial_\theta}{r}(f_k(r)e^{mik\theta})e_r - \partial_r(f_k(r)e^{mik\theta})e_\theta\big)\cdot\big(-\alpha(2-\alpha)\beta r^{-1-\alpha}e_r\big) \biggr).
    \end{align*}
	Simplifying, we get
	\begin{align}
    (L_{ss} \Omega)(r,\theta) 
    &= \sum_{k \in \Z} \biggl( 
	\tfrac{r}{\alpha}\big(f_k''(r)+\tfrac{f_k'(r)}{r}-\tfrac{(mk)^2f_k(r)}{r^2}\big)'+\big(f_k''(r)+\tfrac{f_k'(r)}{r}-\tfrac{(mk)^2f_k(r)}{r^2}\big)\notag \\
	&-mik\beta r^{-\alpha}\big(f_k''(r)+\tfrac{f_k'(r)}{r}-\tfrac{(mk)^2f_k(r)}{r^2}\big)-mikr^{-2-\alpha}\alpha(2-\alpha)\beta f_k(r) \biggr) e^{imk\theta}.\label{HOMOGODE}
	\end{align}
    Hence, the operator $L_{ss}$ leaves the spaces $U_{km}$ invariant and it reduces to the operator $A_{m,k} := (L_{ss})|_{U_{km}}$ given by
	\begin{equation*}
		A_{m,k} (u(r) e^{ikm\theta}) 
		=  \biggl( 
		\frac{r}{\alpha} u'(r)
		+ u(r) - mik\beta r^{-\alpha} u - mik\alpha(2-\alpha)\beta 
		r^{-\alpha-2} f(r)
		\biggr) e^{ikm\theta},
        \end{equation*}
	\begin{equation*} D(A_{k,m}) = \bigl\{ u(r)e^{imk\theta} \in \rL^2_w(\R^2) \colon \tfrac{r}{\alpha} u'
		+ u - mik\beta r^{-\alpha} u - mik\alpha(2-\alpha)\beta 
		r^{-\alpha-2} f
		\in \rL^2(r^{3+\alpha}\mathrm{d}r) \bigr\},
	\end{equation*}
	where $f$ is the solution of 
	\begin{equation}
		u(r) = f''(r) + \frac{f'(r)}{r} - \frac{(mk)^2 f(r)}{r^2} . \label{eq:stream fct}
	\end{equation}
	Using the isometric isomorphism
	\begin{equation*}
		S\colon U_{km} \to \rL^2(\R_+,r^{3+\alpha} \mathrm{d} r) \colon u e^{ikm\theta} \mapsto u
	\end{equation*}
	the operator $B_{m,k} := S A_{k,m} S^{-1}$ on $\rL^2(\R_+,r^{3+\alpha}\mathrm{d}r)$ is given by
	\begin{align*}
		B_{k,m} u &= \frac{r}{\alpha} u'(r)
		+ u(r) - mik\beta r^{-\alpha} u - mik\alpha(2-\alpha)\beta 
		r^{-\alpha-2} f(r), \\
		D(B_{m,k}) &= \bigl\{ U \in \rL^2(\mathbb{R},r^{3+\alpha}\mathrm{d}r) \colon  \tfrac{r}{\alpha} u'
		+ u - mik\beta r^{-\alpha} u - mik\alpha(2-\alpha)\beta f \in \rL^2(\mathbb{R},r^{3+\alpha}\mathrm{d}r) \bigr\}
	\end{align*}
	with $f$ as in \eqref{eq:stream fct}.
	Using the isometric isomorphism
	\begin{equation*}
		T \colon \rL^2(\R_+,r^{3+\alpha}\mathrm{d}r) 
		\to \rL^2(\R,\mathrm{d}s) \colon 
		u \mapsto [ s \mapsto u(e^s)e^{(2+\frac{\alpha}{2})s} ],
	\end{equation*}
	we obtain that the operator $C_{m,k} := T B_{m,k} T^{-1}$ is given by
	\begin{align*}
		(C_{m,k} U)(s) &= \tfrac{1}{\alpha} U'(s) + \bigl(\tfrac{1}{2}-\tfrac{2}{\alpha}\bigr) U(s) - mik\beta e^{-\alpha s} U(s) - mik\alpha(2-\alpha)\beta e^{-\alpha s} \psi(s),\\
		D(C_{m,k}) &= \bigl\{ U \in \rL^2(\mathbb{R}) \colon  \tfrac{1}{\alpha} U' - mik\beta e^{-\alpha \bullet} U - mik\alpha(2-\alpha)\beta e^{-\alpha \bullet} \psi \in \rL^2(\mathbb{R}) \bigr\}
	\end{align*}
	with $\psi(s) := f(e^s)e^{\frac{\alpha}{2}s}$ satisfying
	\begin{equation}
		U(s) = \psi''(s) -\alpha\psi'(s) + \left(\frac{\alpha^2}{4}-(mk)^2\right)\psi(s) . \label{eq:stream fct psi}
	\end{equation}
	Note that, since $S$ and $T$ are isometric isomorphisms, the operator $L_{ss}$ generates a strongly continuous, quasi-contractive semigroup if and only the operators $C_{m,k}$ generate strongly continuous, quasi-contractive semigroups for all $k \in \Z$. Furthermore, the growth bounds satisfy
	\begin{equation*}
		\omega_0(L_{ss}) = \sup_{k \in \Z} \ \omega_0(C_{m,k}) .
	\end{equation*}
	
	\subsection{Dissipativity}
	\label{ssec:dissipative}
	
%		\subsection{Generator property of the operators}
	
	Throughout this section we fix $m \in \N$ and $k \in \Z$. Since we only deal with one operator at the same time we simplify the notation by setting $C := C_{m,k}$. We split 
	\begin{equation*}
		CU = MU+KU + (\tfrac{1}{2}-\tfrac{2}{\alpha}) U, %=: \tilde{M} U + K U .
	\end{equation*}
	%where $\overline{M+K}$ denotes the closure of $M+K$.
	Here $M$ is the main part and $K$ a kernel operator given by
	\begin{align*}
		M U(s)&:= \tfrac{1}{\alpha} U'(s)-mik\beta e^{-\alpha s}U(s), \\
		D(M)&:= \left\{ U \in \mathrm{L}^2(\mathbb{R}) \colon \tfrac{1}{\alpha} U'-mik\beta e^{-\alpha \bullet}U \in \mathrm{L}^2(\mathbb{R}) \right\}, \\
		KU(s) &:= -mik\alpha(2-\alpha)\beta e^{-\alpha s}\psi(s) ,\\
		D(K) &:= \left\{ U \in \mathrm{L}^2(\mathbb{R}) \colon e^{-\alpha \bullet}\psi \in \mathrm{L}^2(\mathbb{R}) \right\} .
	\end{align*}
	First, we deal with the main part $M$.
	\begin{lem}
		\label{lem:sg_mainpart}
		The operator $M$ generates a strongly continuous contraction semigroup $(T_{M}(t))_{t\geq 0}$ on $\rL^2(\mathbb{R})$ given by:
		\begin{equation}
			(T_{M}(t)U)(s) =  e^{mik\beta \bigl(e^{-t-\alpha s} - e^{-\alpha s} \bigr) }U(\alpha^{-1}t + s) .
		\end{equation}
	\end{lem}
	\begin{proof}
		The proof is divided into several steps.
		
		First we show that $T_{M}(t)$ are bounded operators on $\mathrm{L}^2(\mathbb{R})$. To this end, note that 
		\begin{equation*}
			|(T_{M}(t)U)(s)|
			= |U(\alpha^{-1} t+s)| 
		\end{equation*}
		for all $s \in \mathbb{R}$ and $t \geq 0$. Now it follows
		\begin{equation*}
			\| (T_{M}(t)U) \|_{\mathrm{L}^2}^2
			= \int_{\mathbb{R}} |U(\alpha^{-1} t+s)|^2 \mathrm{d} s
			= \int_{\mathbb{R}} |U(x)|^2 \mathrm{d} x 
			= \| U \|_{\mathrm{L}^2}^2 . 
		\end{equation*}
		
		In the next step, we show the semigroup property. Note that $T(0) U = U$ and 
		\begin{align*}
			&\phantom{ = } \ (T_{M}(t_1) T_{M}(t_2)U)(s) \\
			&= T_{M}(t_1) e^{mik\beta (e^{-t_1-\alpha s-e^{-\alpha s})} )}U(\alpha^{-1}t_1 + s) \\
			&= e^{mik\beta (e^{-t_2-\alpha s}-e^{-\alpha s})} e^{mik\beta (e^{-t_1- \alpha(\alpha^{-1} t_2-s)}-e^{-\alpha (\alpha^{-1}t_2+s)})} U(\alpha^{-1}t_1 + \alpha^{-1}t_2 + s) \\
			&= e^{mik\beta (e^{-t_2-\alpha s}-e^{-\alpha s}+ e^{-(t_1+ t_2)-\alpha^{-1}s}-e^{-t_2 - \alpha s})} U(\alpha^{-1}(t_1 +t_2) + s)
			\\
			&= e^{mik\beta (e^{-t_1- t_2-\alpha^{-1}s}-e^{-\alpha s})} U(\alpha^{-1}t_1 + \alpha^{-1}t_2 + s) \\
			&= T_{M}(t_1+t_2)U(s) .
		\end{align*}
		Now, we prove strong continuity of $(T_{M}(t))_{t \geq 0}$. From the first step we know that $(T(t))_{t \geq 0}$ is uniformly bounded on $\mathbb{R}$. For $U \in \mathrm{C}_c(\mathbb{R})$, we obtain
		\begin{align*}
			&\phantom{=} \ \| T_{M}(t)U - U \|_{\mathrm{L}^2}
			\leq C \cdot \| T_{M}(t)U - U \|_{\mathrm{L}^\infty} 
			\\
			&\leq 
			C \sup_{s \in \mathbb{R}}\left| e^{mik\beta \bigl(e^{-t-\alpha s}-e^{-\alpha s} \bigr) }-1 \right| \cdot
			|U(s)|\\
			&+ C \sup_{s \in \mathbb{R}} 
			e^{mik\beta \bigl(e^{-t-\alpha s)}-e^{-\alpha s} \bigr) } \left| U(\alpha^{-1}t+s) - U(s)\right| 
			\\
			&\leq C \sup_{s \in \mathbb{R}}\left| e^{mik\beta \bigl(e^{-t-\alpha s}-e^{-\alpha s} \bigr) }-1 \right| \cdot
			\| U \|_{\mathrm{L}^\infty} \\
			&+ 2 C \sup_{s \in \mathbb{R}} 
			\left| U(\alpha^{-1}t+s) - U(s)\right|,
		\end{align*}
		which converges to $0$ for $t \downarrow 0$. Hence, we obtain 
		\begin{equation*}
			\| T_{M}(t)U - U \|_{\mathrm{L}^2} \to 0 
		\end{equation*}
		for $t \downarrow 0$ for all $U \in \mathrm{C}_c(\mathbb{R})$ and by density for all $U \in \mathrm{L}^2(\mathbb{R})$.
		Hence $(T_{M}(t))_{t \geq 0}$ is a strongly continuous semigroup on $\mathrm{L}^2(\mathbb{R})$.
		
		\medskip 
		
		Finally, we show that $M$ is the generator of $(T_{M}(t))_{t \geq 0}$ on $\mathrm{L}^2(\mathbb{R})$.
		Product rule implies
		\begin{align*}
			&\frac{d}{dt}\bigg|_{t = 0} T_{M}(t)U(s) \\
			= \ &e^{mik\beta(e^{-t-\alpha s}-e^{-\alpha s}) } \alpha^{-1} U'(\alpha^{-1}t + s)\bigg|_{t = 0}
			- mik\beta e^{-t-\alpha s} e^{mik\beta(e^{-t-\alpha s}-e^{-\alpha s})} U(\alpha^{-1}t+s)
			\bigg|_{t = 0} \\
			= \ &\alpha^{-1} U'(s) - mik\beta e^{-\alpha s}U(s)
		\end{align*}
		for $U \in D(M)$.
	\end{proof}
	
	\begin{lem}\label{lem:domain M}
		The domain of $M$ is given by
		\begin{equation*}
			D(M) = \left\{ e^{mik\beta e^{-\alpha \bullet}}U \in \mathrm{H}^1(\R)  \right\} . 
		\end{equation*}
	\end{lem}
	\begin{proof}
		Note that
		\begin{align*}
			(e^{mik\beta e^{-\alpha \bullet}}U)'(s)
			&= e^{mik\beta e^{-\alpha s}}(U'(s) 
			- mik\alpha \beta e^{-\alpha s}U(s)) \\
			&= \alpha e^{mik\beta e^{-\alpha s}}(\tfrac{1}{\alpha} U'(s) 
			- mik \beta e^{-\alpha s}U(s)) 
		\end{align*}
		and hence
		\begin{align*}
			\left|e^{-mik\beta e^{-\alpha s}}U(s)\right|^2 &= 
			|U(s)|^2, \\
			\left| (e^{-mik\beta e^{-\alpha \bullet}}U)'(s)\right|^2
			&= \left| \alpha e^{-mik\beta e^{-\alpha s}}(\frac{1}{\alpha} U'(s) 
			- mik \beta e^{-\alpha s}U(s)) \right|^2
			\\
			&= \alpha^2 \cdot  \left|\frac{1}{\alpha} U'(s) 
			- mik \beta e^{-\alpha s}U(s) \right|^2 .
		\end{align*}
		This implies $U \in \rL^2(\mathbb{R})$ if and only if 
		$e^{-mik\beta e^{-\alpha \bullet}}U \in \rL^2(\R)$ and 
		$\tfrac{1}{\alpha} U'-mik\beta e^{-\alpha \bullet}U \in \rL^2(\R)$ if and only if $(e^{-mik\beta e^{-\alpha \bullet}}U)' \in \rL^2(\R)$.
	\end{proof}
	
	%Next we add the "drift" term. In particular, we consider the operator 
	%\begin{equation*}
	%    \tilde{M} U(s) = \tfrac{1}{\alpha} U'(s)+(1-\tfrac{1}{\alpha})U(s)-mik\beta e^{-\alpha s}U(s).    
	%\end{equation*}
	%
	%\begin{lem}
	%\label{EXPLICITSEMI2}
	%    The operator $\tilde{A}_1$ generates a strongly continuous quasi-contractive semigroup $(T_{\tilde{M}}(t))_{t\geq 0}$ on $\rL^2(\mathbb{R})$ with growth bound $\omega_0(\tilde{M}) = 1-\frac{1}{\alpha}< 0$.
	%\end{lem}
	%\begin{proof}
	%    The semigroup $(T_{\tilde{M}}(t))_{t \geq 0}$ is obtained 
	%    by rescaling
	%    \begin{equation*}
		%       T_{\tilde{M}}(t)
		%        = e^{(1-\frac{1}{\alpha})t}T_{M}(t) 
		%    \end{equation*}
	%    and the result follows by \cite[Chapter~II.2.2]{EN}.
	%\end{proof}

	It remains to deal with the operator 
	\begin{equation*}
		KU(s) = -mik\alpha(2-\alpha)\beta e^{-\alpha s}\psi(s) , \qquad  D(K) := \left\{ U \in \rL^2(\R) \colon e^{-\alpha \bullet}\psi\in \rL^2(\R) \right\} .
	\end{equation*}
	
	We recall the following criterion for dissipativity on Hilbert spaces from \cite[Proposition II.3.23]{EN}:
	\begin{lem}
		\label{DISSDEF}
		An operator $A \colon D(A) \subset H \to H$ on a Hilbert space $H$ is dissipative if and only if 
		\begin{equation*}
			\mathrm{Re}\langle Ax, x \rangle \leq 0
		\end{equation*}
		for all $x\in D(A)$.
	\end{lem}
	
	The choice of the weighted $\rL^2_w$ space guarantees the following result.
	%We review some material from
	%Consider the following real-valued kernel:
	%$$K_1(t,s):= \bigg( e^{-(km+1+\alpha)(t-s)}\chi_{(0,\infty)}(t-s) + e^{(km-(1+\alpha))(t-s)}\chi_{(-\infty, 0)}(t-s)\bigg).$$
	%We define $\Phi_1$ to be the bounded operator on $\rL^2$ given by integration against $K_1$. In particular, $\Phi_1$ is defined as
	%$$\Phi_1(U)(t) := -\frac{1}{2km}\int_{\mathbb{R}} K_1(t,s)U(s)ds.$$
	%
	%Consider the following operator:
	%$$\tilde{B}U(s) := -mik\alpha(2-\alpha)\beta e^{-\alpha s}\Phi_1(U)(s).$$
	%We shall now prove that:
	\begin{lem}
		\label{BDISS}
        We have $\Re \langle KU, U \rangle = 0$ for $U \in D(K)$. In particular, 
		the operator $K$ is dissipative on $\rL^2(\mathbb{R})$.
	\end{lem}
	\begin{proof}
		We proceed directly using Equation \eqref{eq:stream fct psi}:
        \begin{align*}
        \langle KU, U\rangle &= \int_{\mathbb{R}} -mik\alpha(2-\alpha)\beta e^{-\alpha s} \psi(s) \overline{U(s)}ds \\
        &= -mik\alpha(2-\alpha)\beta\int_{\mathbb{R}} e^{-\alpha s}\psi(s)\left(\overline{\psi''(s)}-\alpha\overline{\psi'(s)}+\left(\tfrac{\alpha^2}{4}-(mk)^2\right)\overline{\psi(s)}\right)ds.
        \end{align*}
        Therefore, after an integration by parts:
        \begin{equation*}
        \Re \langle KU, U\rangle = 0.
        \end{equation*}
        \qedhere
	\end{proof}
	Let $A:D(A)\subset X\to X$ be a densely defined linear operator on a Banach space $X$. We call a densely defined operator $B:D(B)\subset X\to X$ {\textbf{relatively $A$-bounded}} if and only if $D(A)\subset D(B)$ and if there exists constants $a,b>0$ such that 
	\begin{equation}\label{RELBOUND}\|Bx\| \leq a\|Ax\| + b\|x\|.\end{equation}
	In this case, we call the number
	$$a:=\inf_{a'\geq 0} \{ \textrm{there exists } b\geq 0 : \textrm{Equation } \eqref{RELBOUND} \textrm{ holds}\}$$
	the {\textbf{$A$-bound of $B$}}. 
	\begin{lem}\label{lem:domain K}
		$K$ is relatively bounded with respect to $M$ with relative bound 
        \begin{equation}\label{RELBOUNDVALUE}
        a = \frac{2\alpha(2-\alpha)}{(2m^2-\alpha^2)}.
        \end{equation}
        In particular, we have $D(M) \subseteq D(K)$. 
	\end{lem}
	\begin{proof}
		We set $\tilde{U}(s) := e^{mik\beta e^{-\alpha s}}U(s)$ and $\tilde{\psi}(s) := -mik\alpha(2-\alpha)\beta e^{-a s} \psi=KU$.
		From \eqref{eq:stream fct psi} we obtain 
	\begin{align}
			\tilde{\psi}''(s) + \alpha \tilde{\psi}'(s)
			+ (\tfrac{\alpha^2}{4}-(km)^2) \tilde{\psi}(s) 
			&= -mik\alpha(2-\alpha)\beta e^{-\alpha s} U(s) \
			\notag \\
            &= -mik\alpha(2-\alpha)\beta e^{-\alpha s} e^{-mik\beta e^{-\alpha s}}\tilde{U}(s)  \label{eq:psi tilde}\\
            &= -(2-\alpha)\left(e^{mik\beta e^{-\alpha s}}\right)' \tilde{U}(s) \notag .
		\end{align}
		Multiplying \eqref{eq:psi tilde} by $-\overline{\tilde{\psi}}$ and integrating by parts we obtain
    $$\|\tilde{\psi}'\|_{\rL^2}^2+ ((km)^2-\tfrac{\alpha^2}{4})\|\tilde{\psi}\|_{\rL^2}^2 = \alpha \int \tilde{\psi}'\overline{\psi}ds -(2-\alpha)\int \left(e^{mik\beta e^{-\alpha s}}\right)' \tilde{U}(s)\overline{\tilde{\psi}(s)}ds.$$
We use the Cauchy-Schwarz inequality for the first term on the right-hand side and integrate by parts once more for the second term. This yields
\begin{align*}
&\|\tilde{\psi}'\|_{\rL^2}^2+ ((km)^2-\tfrac{\alpha^2}{4})\|\tilde{\psi}\|_{\rL^2}^2 \\
\leq \, &\alpha \|\tilde{\psi}'\|_{\rL^2}\|\tilde{\psi}\|_{\rL^2} +(2-\alpha)\left|\int e^{mik\alpha\beta e^{-\alpha s}} \tilde{U}'(s)\overline{\tilde{\psi}(s)}ds\right|+(2-\alpha)\left|\int e^{mik\alpha\beta e^{-\alpha s}} \tilde{U}(s)\overline{\tilde{\psi}'(s)}ds\right|.
\end{align*}
Therefore, by applying Cauchy-Schwarz once more we have
\begin{align*}
&\|\tilde{\psi}'\|_{\rL^2}^2+ ((km)^2-\tfrac{\alpha^2}{4})\|\tilde{\psi}\|_{\rL^2}^2 \\
\leq \, &\alpha \|\tilde{\psi}'\|_{\rL^2}\|\tilde{\psi}\|_{\rL^2} +(2-\alpha) \|\tilde{U}'\|_{\rL^2}\|\tilde{\psi}\|_{\rL^2}+(2-\alpha)\|\tilde{U}\|_{\rL^2}\|\tilde{\psi}'\|_{\rL^2}.
\end{align*}
We aim for best possible constant for the relative bound, so let us apply the Peter-Paul inequality with some $\epsilon_1, \epsilon_2,\epsilon_3>0$:
\begin{align*}
&\|\tilde{\psi}'\|_{\rL^2}^2+ ((km)^2-\tfrac{\alpha^2}{4})\|\tilde{\psi}\|_{\rL^2}^2 \\
\leq &\alpha\frac{\epsilon_1^2}{2} \|\tilde{\psi}'\|_{\rL^2}^2+\alpha\frac{1}{2\epsilon_1^2}\|\tilde{\psi}\|_{\rL^2}^2 +(2-\alpha)\frac{\epsilon_2^2}{2}\|\tilde{U}'\|_{\rL^2}^2+(2-\alpha)\frac{1}{2\epsilon_2^2}\|\tilde{\psi}\|_{\rL^2}^2\\
&+(2-\alpha)\frac{1}{2\epsilon_3^2}\|\tilde{U}\|_{\rL^2}^2+(2-\alpha)\frac{\epsilon_3^2}{2}\|\tilde{\psi}'\|_{\rL^2}^2.
\end{align*}
We conclude:
\begin{align*}
\left(1-\tfrac{\alpha\epsilon_1^2}{2}-\tfrac{(2-\alpha)\epsilon_3^2}{2}\right)\|\tilde{\psi}'\|_{\rL^2}^2+\left((km)^2-\tfrac{\alpha^2}{4} -\tfrac{\alpha}{2\epsilon_1^2}-\frac{2-\alpha}{2\epsilon_2^2}\right)\|\tilde{\psi}\|_{\rL^2}^2
\\\leq 
\frac{(2-\alpha)\epsilon_2^2}{2}\|\tilde{U}'\|_{\rL^2}^2+\frac{2-\alpha}{2\epsilon_3^2}\|\tilde{U}\|_{\rL^2}^2.
\end{align*}
The optimal choice of $\epsilon_1, \epsilon_2, \epsilon_3$ gets us for any $\epsilon>0$ something like:
$$\|\tilde{\psi}\|_{\rL^2}^2 \leq \frac{4(2-\alpha)^2(1+\epsilon)}{(2m^2-\alpha^2)^2}\|\tilde{U}'\|_{\rL^2}^2 + \frac{C_{\alpha,m}}{\epsilon}\|\tilde{U}\|_{\rL^2}^2.$$
We require the conditions $\alpha \in (0,2)$ since $m\geq 2$ and for $m=1$ we require $\alpha \in (0,\sqrt{2})$.
From this we get
$$\|KU\|_{\rL^2}= \|\tilde{\psi}\|_{\rL^2} \leq \frac{2\alpha(2-\alpha)(1+\epsilon)}{(2m^2-\alpha^2)}\|MU\|_{\rL^2}+\frac{C_{\alpha,m}}{\epsilon}\|\tilde{U}\|_{\rL^2}.$$
So the optimal relative bound we get is
$$a = \frac{2\alpha(2-\alpha)}{(2m^2-\alpha^2)} < 1,$$
for all $\alpha \in (0,2)$ if $m\geq 2$, and is $<1$ for all $\alpha \in (0,2-\sqrt{2})$ if $m=1$.
	\end{proof}
	
	Now, we are in the position to define 
	the sum of $M$ and $K$ by
	\begin{equation*}
		C_0 U := MU + KU, \qquad D(C_0) = D(M) .
	\end{equation*}

    	\subsection{Conclusion}
	\label{sec:conclusion}
		
	We recall the following standard perturbation theorem (cf. \cite[Theorem III.2.7]{EN}):
	\begin{thm}
		\label{LP}
		Let $(A,D(A))$ be the generator of a contraction semigroup and assume $(B,D(B))$ to be dissipative and $A$-bounded with $A$-bound $a_0<1$. Then $(A+B, D(A))$ generates a contraction semigroup.
	\end{thm}
	
	Now, we are ready to conclude our main result.
	Combining the perturbation theorem from \autoref{LP} with the relative bound we found in \eqref{RELBOUNDVALUE}, we conclude the following result.
	
	\begin{cor}\label{prop:C_0 sg}
		The operator $C_0$ generates a strongly continuous contractive semigroup on $\rL^2(\R)$ that we denote $(T_{C_0}(t))_{t\geq 0}$.
	\end{cor}
	
	Finally, we add the drift term $\tfrac{1}{2}-\tfrac{2}{\alpha}$ to $C_0$ and obtain the following result.
	
	\begin{cor}\label{cor:C_mk sg}
		The operator $C$ generates a strongly continuous, quasi-contractive semigroup on $\rL^2(\mathbb{R})$ with growth bound $\omega_0(C) \leq \tfrac{1}{2}-\tfrac{2}{\alpha}$.
	\end{cor}
	\begin{proof}
		Recall that 
		\begin{equation*}
			C U = C_0 + (\tfrac{1}{2}-\tfrac{2}{\alpha})U .
		\end{equation*}
		Using \autoref{prop:C_0 sg} we obtain by \cite[Section II.2.2]{EN} that $C$ generates the rescaled semigroup
		\begin{equation*}
			T_C(t) = e^{\bigl(\frac{1}{2}-\frac{2}{\alpha}\bigr)t} \, T_{C_0}(t) 
		\end{equation*}
		and the claim is an immediate consequence. 
	\end{proof}
	
	%\subsection{Proof of \autoref{thm:mainss}}
	\begin{proof}[Proof of \autoref{thm:mainss}]
	We have seen in \autoref{ssec:prelims} that it suffices to prove that all $C_{m,k}$ generate strongly continuous, quasi-contractive semigroups; this was proven in \autoref{cor:C_mk sg}. Moreover their growth bounds are satisfy $\omega_0(C_{m,k}) = \tfrac{1}{2}-\tfrac{2}{\alpha}$. Hence, we obtain 
	\begin{equation*}
		\omega_0(L_{ss}) = \sup_{k \in \Z} \ \omega_0(C_{m,k}) \leq \tfrac{1}{2}-\tfrac{2}{\alpha} .
	\end{equation*}
	Using \autoref{cor:spectral bound1} in conjunction with \cite[Proposition IV.2.2.]{EN} we conclude
	\begin{equation*}
		\tfrac{1}{2} - \tfrac{2}{\alpha} \leq s(L_{ss}) \leq \omega_0(L_{ss}) \leq \tfrac{1}{2}-\tfrac{2}{\alpha}
	\end{equation*}
	and the claim is proven. 
	\end{proof}

	\subsection{Lower Bound for the Spectral Bound}
	\label{ssec:spectral bound1}
	
	In this subsection we show that  the line $\textrm{Re}(\lambda)=1/2-2/\alpha$ is contained in the spectrum of $L_{ss}$ when defined on $\rL^2_{w,m}$. More precisely, we have the following 
	
	\begin{prop}
		$\lambda = \frac{1}{2}-\frac{2}{\alpha}-i \gamma\in \sigma_{\alpha}(L_{ss},\rL^2_{w,m})$ for all $\alpha\in (0,2)$, $m\geq 2$ or $\alpha \in (0,2-\sqrt{2})$ and $m=1$, and $\gamma\in \mathbb{R}$. 
	\end{prop}
    \begin{proof}  
		From \autoref{ssec:prelims} we know that 
		it suffices to show that $\frac{1}{2}-\frac{2}{\alpha}+i\gamma\in\sigma_{\alpha}(L_{ss},U_{km})$ for some $k\in \mathbb{Z}$. We do this for $k=0$, in which case the resolvent equation for $U(t)$, as above, simplifies to:
		\begin{equation}\label{DERIVBAD1}
			U'(t)+i\alpha\gamma U(t)=\alpha G(t).
		\end{equation}
        Integrating, Equation \eqref{DERIVBAD1} becomes the integral equation
		$$U(t) = ce^{-i\alpha\gamma t} + \alpha e^{-i\alpha\gamma t}\int_0^t e^{i\alpha \gamma s}G(s) ds,$$
		where $c$ is a constant that must be chosen so that $U\in \rL^2(\mathbb{R})$. Consider the function
		$$G(t) := \frac{e^{-i\alpha\gamma s}}{\abs{t}+1}\in \rL^2(\mathbb{R}).$$
		It is not difficult to see that 
		$$ce^{-i\alpha\gamma t} + \alpha e^{-i\alpha\gamma t}\int_0^t e^{i\alpha \gamma s}G(s) ds = e^{-i\alpha\gamma t}\big(c + \alpha\cdot\textrm{sgn}(t)\log(1+\abs{t})\big),$$
		is not in $\rL^2(\mathbb{R})$ for any choice of the constant $c$. This implies that $\frac{1}{2}-\frac{2}{\alpha}-i\gamma \in \sigma_\alpha(L_{ss}, U_0)$ for all $\gamma\in \mathbb{R}$.
	\end{proof}
	
	As a direct consequence, we obtain for the spectral bound of $L_{ss}$ the following result.
	
	\begin{cor}\label{cor:spectral bound1}
		The spectral bound of $L_{ss}$ is bounded from below by
		\begin{equation*}
			\frac{1}{2}-\frac{2}{\alpha} \leq s(L_{ss}) .
		\end{equation*}
	\end{cor}

\section{Linear Stability of the Power-Law Vortex in Physical Coordinates}
	\label{sec:linear stability}
In this section we deal with the linearization of the Euler equation around the power-law vortex in physical coordinates.
Our goal is to prove the following result.
\begin{thm}
    \label{thm:mainp}
    Let $\beta\in \mathbb{R}\setminus\{0\}$ and $m\geq 2$ and $\alpha \in (0,2)$ or $m=1$ and $\alpha \in (0,2-\sqrt{2})$. 
    Let $\overline{\Omega}(\xi) = \beta(2-\alpha) \abs{\xi}^{-\alpha}$ be the radial power-law vortex that solves the unforced Euler equations in self-similar coordinates and let $L_{\phi}$ be the linearization of the Euler equations around $\overline{\Omega}$. Then the operator $L_{\phi}$ generates a continuous group of isometries on $\rL^2_{w,m}(\R^2)$ with
    $$\sigma(L_\phi) \subset \{z \in \mathbb{C} : \textrm{Re}(z)=0\}.$$
\end{thm}

    As remarked earlier, the background vortex $\overline{\Omega}(x)= \beta\abs{x}^{-\alpha}$ is a stationary solution of the Euler equations in physical coordinates as well. The linearization of the Euler equations in vorticity form around $\overline{\Omega}$ is precisely
	\begin{equation}
		\label{LINEAREULERP}
		\partial_t \Omega +\overline{V}\cdot\nabla \Omega + V\cdot\nabla \overline{\Omega}=0\end{equation}
	$$K_2\ast V = \Omega.$$
	Alternately, we can write the linear Equation \eqref{LINEAREULERP} as
	$$(\partial_t - L_\phi)\Omega =0,$$
	where the operator $L_\phi \colon D(L_\phi) \subset \rL^2(\R^2,|x|^{2+\alpha} \mathrm{d} x) \to \rL^2(\R^2,|x|^{2+\alpha} \mathrm{d} x)$ given by
	\begin{align*}
		L_\phi \Omega 
		&= -\overline{V}\cdot\nabla \Omega - V\cdot\nabla \overline{\Omega}, \\
		D(L_\phi)
		&= \{ \Omega \in \rL^2(\R^2,|x|^{2+\alpha} \mathrm{d} x) \colon L_\phi \Omega \in \rL^2(\R^2,|x|^{2+\alpha} \mathrm{d} x) \} .
	\end{align*}
    Following the arguments from \autoref{ssec:prelims} we see that it suffices to study the operators
\begin{align*}
		(G_{m,k} U)(s) &=  - mik\beta e^{-\alpha s} U(s) - mik\alpha(2-\alpha)\beta e^{-\alpha s} \psi(s),\\
		D(G_{m,k}) &= \bigl\{ U \in \rL^2(\mathbb{R}) \colon e^{-\alpha \bullet} U +\alpha(2-\alpha) e^{-\alpha \bullet} \psi \in \rL^2(\mathbb{R}) \bigr\}
\end{align*}
	with $\psi(s)$ satisfying
	\begin{equation*}
		U(s) = \psi''(s) -\alpha\psi'(s) + (\tfrac{\alpha^2}{4}- (mk)^2) \psi(s) . 
	\end{equation*}
    In order to simplify the notation we denote $G_{m,k}$ by $G$ whenever it is clear from the context which $m$ and $k$ are meant. 
    Again we split the operator
    \begin{equation*}
        G = \tilde{M} + K
    \end{equation*}
    in a main part $\tilde{M}$ given by
    \begin{equation*}
    \tilde{M}U(s) := -mik\beta e^{-\alpha s} U(s),
    \qquad D(\tilde{M}) = D(G) 
    \end{equation*}
    and a perturbation part $K$ given by
    \begin{equation*}
		KU(s) = -mik\alpha(2-\alpha)\beta e^{-\alpha s}\psi(s) , \qquad  D(K) := \left\{ U \in \rL^2(\R) \colon e^{-\alpha \bullet}\psi\in \rL^2(\R) \right\} .
    \end{equation*}
    Note that the operator $K$ is the same operator as in \autoref{sec:linear stability self-similar}.
    %Next, we denote
%	$\tilde{L}_{p,k}$ to be the linearization of the Euler equations around $\overline{\Omega}$ in physical coordinates, after performing a separation of variables and the transform $r\to e^t$. Explicitly, it is given by
%	$$\tilde{L}_{p,k}U(s) := \tilde{M}U(s)+\tilde{B}U(s).$$
%

	\begin{lem}\label{lem:relativebound}
    Suppose $m\geq 2$ and $\alpha \in (0,2)$ or $m=1$ and $\alpha \in (0,2-\sqrt{2})$. 
    Then the operator $K$ is relatively $\tilde{M}$-bounded with $\tilde{M}$-bound strictly less than one.
	\end{lem}
	\begin{proof}
		It suffices to prove that 
		\begin{equation}\label{RELBOUNDGOAL}\alpha(2-\alpha)\|e^{-\alpha s}\psi(s)\|_{\rL^2}\leq a\|e^{-\alpha s}U(s)\|_{\rL^2} + b\|U(s)\|_{\rL^2}\end{equation}
		for some constants $0\leq a<1$ and $b\geq 0$. Let $\tilde{\psi} = \alpha(2-\alpha)e^{-\alpha s}\psi$ and $\tilde{U}(s)=e^{-\alpha s}U(s)$. The ordinary differential equation satisfied by $\tilde{\psi}$ is
        $$\tilde{\psi}''+\alpha\tilde{\psi}'+\left(\tfrac{\alpha^2}{4}-(km)^2\right)\tilde{\psi}=(2-\alpha)\alpha e^{-\alpha s}U(s)=(2-\alpha)\alpha \tilde{U}. $$
        We test the equation against $\overline{\tilde{\psi}}$, integrate by parts, and use the Cauchy-Schwarz inequality to get:
$$\|\tilde{\psi}'\|_{\rL^2}^2 + \left((km)^2-\tfrac{\alpha^2}{4}\right)\|\tilde{\psi}\|_{\rL^2}^2 \leq \alpha\|\tilde{\psi}'\|_{\rL^2}\|\tilde{\psi}\|_{\rL^2} + \alpha(2-\alpha)\|\tilde{U}\|_{\rL^2}\|\tilde{\psi}\|_{\rL^2}.$$
Once again, we use the Peter-Paul inequality to get the optimal relative bound. To that end, let $\epsilon_1,\epsilon_2>0$. We get:
$$\|\tilde{\psi}'\|_{\rL^2}^2 + \left((km)^2-\tfrac{\alpha^2}{4}\right)\|\tilde{\psi}\|_{\rL^2}^2 \leq \frac{\alpha\epsilon_1}{2}\|\tilde{\psi}'\|_{\rL^2}^2+\frac{\alpha}{2\epsilon_1}\|\tilde{\psi}\|_{\rL^2}^2 + \frac{\alpha(2-\alpha)\epsilon_2}{2}\|\tilde{U}\|_{\rL^2}^2+\frac{\alpha(2-\alpha)}{2\epsilon_2}\|\tilde{\psi}\|_{\rL^2}^2.$$
Whence
$$\left(1-\frac{\alpha\epsilon_1}{2}\right)\|\tilde{\psi}'\|_{\rL^2}^2 + \left((km)^2-\frac{\alpha^2}{4}-\frac{\alpha}{2\epsilon_1}-\frac{\alpha(2-\alpha)}{2\epsilon_2}\right)\|\tilde{\psi}\|_{\rL^2}^2 \leq \frac{\alpha(2-\alpha)\epsilon_2}{2}\|\tilde{U}\|_{\rL^2}^2.$$
The choice of optimal constant gets us:
$$\|\tilde{\psi}\|_{\rL^2} \leq \frac{2 (2-\alpha) \alpha }{2m^2-\alpha^2}\|\tilde{U}\|_{\rL^2},$$
so long as $m\geq 2$ and $\alpha \in (0,2)$ or $m=1$ and $\alpha \in (0,\sqrt{2})$. 
As before, we know that the relative bound is $a<1$ if $m\geq 2$ or if $m=1$ and $\alpha \in (0,2-\sqrt{2})$.
	\end{proof}
	
	\begin{lem}
		\label{lem:mcontract}
        The operators $\pm\tilde{M}$ are the generators of contraction semigroups on $\rL^2(\R)$.
	\end{lem}
	\begin{proof}
		Observe that $\pm \tilde{M}$ are multiplication operators, given by multiplication by purely imaginary functions. Hence,
		$$\sigma(\pm\tilde{M})= \{z\in\mathbb{C} : \textrm{Re}(z)=0\}.$$
		Now, \cite[Section II.2.9]{EN} implies that $\pm\tilde{M}$ generate strongy continuous contraction semigroups on $\rL^2(\R)$.
	\end{proof}

        We need in the sequel the following generation theorem for strongly continuous groups of isometries, cf. \cite[Corollary II.3.7 \& Proposition II.3.11 (Generation Theorem for Groups)]{EN}:

        \begin{thm}\label{thm:generation groups}
            For a linear operator $(A,D(A))$ on a Banach space $X$ the following statements are equivalent:
            \begin{enumerate}
                \item[(a)] the operator $A$ generates a strongly continuous group of isometries on $X$.
                \item[(b)] the operators $\pm A$ generate strongly continuous contraction semigroups on $X$.
            \end{enumerate}
        \end{thm}

        Now, we conclude
    
	\begin{prop}
		\label{prop:CONTGROUP}
		Suppose $\alpha \in (0,2)$ and $m\geq 2$ or $\alpha \in (0,2-\sqrt{2})$ and $m=1$. Then the operator $G$ generates a strongly continuous group of isometries $(T_{G}(t))_{t\geq 0}$ on $\rL^2_{w,m}(\mathbb{R})$.
	\end{prop}
	\begin{proof}
        By \autoref{thm:generation groups} we have to show that $\pm G$ 
        generate strongly continuous contraction semigroups on $\rL^2(\R)$.
        Note that 
        \begin{equation*}
            \pm G
            = \pm (\tilde{M} + K) . 
        \end{equation*}
        By \autoref{lem:mcontract}, we know that $\pm\tilde{M}$ are generators of contraction semigroups. By \autoref{lem:relativebound}, we know that $K$ is relatively bounded with respect to $\tilde{M}$ with $\tilde{M}$-bound strictly less than one. 
        Note that relative boundedness does not depend on the sign and hence $-K$ is relatively $(-\tilde{M})$-bounded with the same bound.
        Moreover, by \autoref{BDISS}, we know that $K$ is a dissipative operator. Further, we see from \autoref{BDISS} that 
        \begin{equation*} 
        \Re \langle - K U,U\rangle = - \Re \langle K U,U \rangle = 0 
        \end{equation*}
        and by \autoref{DISSDEF} $-K$ is dissipative operator, too. 
        Therefore, by \autoref{LP}, we have that $\tilde{M}+K$ and $-\tilde{M}-K$ are generators of contraction semigroups on $\rL^2(\R)$, which proves the claim. 
        %which we may denote by $(T_{\tilde{L}_{p,k}}(t))_{t\geq 0}$.
		%
		%However, $(T_{\tilde{L}_{p,k}}(t))_{t\geq 0}$ actually generates a continuous contraction group, due to the reversibility of the problem on the PDE level. Thus, since it is straightforward to see that the only continuous contraction groups are continuous groups of isometries, we conclude that $\tilde{L}_{p,k}$ generates a strongly continuous group of isometries. 
	\end{proof}
	
	    We recall the weak spectral mapping theorem from \cite[Theorem IV.3.16]{EN}:
	\begin{thm}[Weak Spectral Mapping]
		\label{thm:WSM}
		Let $(T(t))_{t\in \mathbb{R}}$ be a bounded strongly continuous group on a Banach space $X$ with generator $A$. Then the following holds:
		$$\sigma(T(t)) = \overline{e^{t\sigma(A)}} \quad \forall t\in \mathbb{R}.$$
	\end{thm}
	
	\begin{lem}
		\label{lem:CTS}
		Let $(T(t))$ be a continuous group of isometries. Then
		$$\sigma(T(t))\subset \mathbb{S}^1 \subset \mathbb{C}.$$
	\end{lem}
	\begin{proof}
		Since each operator $T(t)$ is an isometry, we certainly have $\abs{\lambda} \leq \|T(t)\|\leq 1$ for every $\lambda \in \sigma(T(t))$. Since $T(t)$ is a group, $T(t)$ is invertible for every $t$ and $\lambda^{-1}\in \sigma(T(t)^{-1})$ for each $t$. Since $T(t)^{-1}=T(-t)$ is an isometry, we have
		$$\abs{\lambda}^{-1}\leq \|T(t)^{-1}\| \leq 1$$
		whence $\abs{\lambda}\geq 1$. We conclude that $\abs{\lambda}=1$ for every $\lambda \in \sigma(T(t))$.
	\end{proof}
	
	\begin{proof}[Proof of \autoref{thm:mainp}]
		First we show that for all $k\in\mathbb{Z}$ we have
		$$\sigma(G_{k,m})\subset \{z\in \mathbb{C} : \textrm{Re}(z)=0\}.$$
		
		Recall that in \autoref{prop:CONTGROUP}, we proved that $\tilde{L}_{p,k}$ generates a strongly continuous group of isometries on $\rL^2(\mathbb{R})$ that we denoted $(T_{G_{p,k}}(t))_{t\geq 0}$. 
		By \autoref{lem:CTS}, we know that $\sigma(T_{G_{p,k}}(t))\subset \mathbb{S}^1$. By the Weak Spectral Mapping Theorem, or \autoref{thm:WSM}, we have
		$$e^{t\sigma(G_{p,k})}\subset \overline{e^{t\sigma(G_{p,k})}}= \sigma(T_{G_{p,k}}(t))\subseteq \mathbb{S}^1.$$
		However, the only complex numbers that map into the unit circle by the exponential map lie on the imaginary axis. Since $t\in\mathbb{R}$, we conclude that
		$$\sigma(G_{p,k}) \subset \{z\in\mathbb{C} : \textrm{Re}(z)=0\}$$
		for all $k\in \mathbb{Z}$. Now, as before, this implies the conclusion of the desired theorem.
	\end{proof}
    
	\appendix

    \section{Spectral Bound for the Linearization in Lebesgue Spaces}
\label{sec:Surjection}

In the case when the underlying space of the analysis is the un-weighted $\rL^2$ space $\rL^2_m$, the nonlocal part of the operator, $K$, is no longer symmetric. Therefore, we proceed differently, using a fixed point argument to show that the resolvent is surjective. We remark that, as we prove in the appendix, the point spectrum is always empty (with no symmetry condition required).

Recall that we wrote the linearization of \eqref{SEULER} around $\overline{\Omega}$ as 
	$$(\partial_\tau -L_{ss})\Omega=0$$
	where the linear operator $L_{ss}$ is given by:
	$$L_{ss} \Omega= \Omega +\big(\frac{\xi}{\alpha}-\overline{V}\big)\cdot\nabla \Omega - \nabla^\perp \Delta^{-1}\Omega \cdot \nabla\overline{\Omega}.$$
    The goal of this section is to prove \autoref{thm:MAINMAINBANACH} and \autoref{BUSYPROP}.
	
	Henceforth we let $g$ be an arbitrary function in $\rL^2_m$ when $q=2$ or $U_{k,q}$ when $q\neq 2$.
	Proving that $L_{ss}-\lambda I$ is surjective with bounded resolvent is equivalent to proving that the operator $R_\lambda(g):= (L_{ss}-\lambda I)^{-1}(g)$ exists and is bounded as an operator into $D(L_{ss})$. We shall succeed in demonstrating surjectivity for $\lambda$ sufficiently large by using a fixed point argument in a class of sufficiently symmetry functions. Proving that $L_{ss}-\lambda I$ is injective amounts to proving that $L_{ss}-\lambda I= g$ has at most one solution. The uniqueness of the solution (or the injectivity of $L_{ss}-\lambda I$) follows from the uniqueness of the fixed point. We shall now reduce our problem to showing that the unique solution of an integral equation is in some Lebesgue space with quantified $\rL^q$ bound. 
	
	In particular, that $L_{ss}-\lambda I$ is surjective is equivalent to stating that the inhomogeneous ordinary differential equation
	$$
	\frac{r}{\alpha}\big(f_k''(r)+\frac{f_k'(r)}{r}-\frac{(mk)^2f_k(r)}{r^2}\big)'+(1-(\lambda_1+i\lambda_2))\big(f_k''(r)+\frac{f_k'(r)}{r}-\frac{(mk)^2f_k(r)}{r^2}\big)$$
	\begin{equation}
		\label{INHOMOG}-mik\beta r^{-\alpha}\big(f_k''(r)+\frac{f_k'(r)}{r}-\frac{(mk)^2f_k(r)}{r^2}\big)-mikr^{-2-\alpha}\alpha(2-\alpha)\beta f_k(r)=g_k(r)
	\end{equation}
	has a solution $f_k$ with $f_k''(r)+\frac{f_k'(r)}{r}-\frac{(mk)^2f_k(r)}{r^2}\in \rL^q(\mathbb{R}, rdr)$ for any given function $g_k\in \rL^q(\mathbb{R},rdr)$. That $L_{ss}-\lambda I$ is injective is equivalent to stating that the homogeneous ordinary differential equation Equation \eqref{HOMOGODE} has no non-trivial solution or, equivalently, that $L_{ss}-\lambda I =g$ has a unique solution for every choice of the function $g$.
	
	Henceforth, we drop the subscript notation showing dependence on $k$ and simply write $f(r)$ and $g(r)$. Moreover, we define:
	$$u(r):= f''(r) +\frac{f'(r)}{r} -\frac{(mk)^2 f(r)}{r^2},$$
	which is the quantity we wish to control in $\rL^q(rdr)$. We now perform the change of variables $e^t =r$ and define:
	$$\psi(t):= f(e^t)e^{(2/q-2)t} \quad G(t):= g(e^t)e^{2t/q} \quad U(t):= u(e^t)e^{2t/q}.$$
	These functions are chosen so that 
	$$g(r) \in \rL^q(rdr) \quad \Longleftrightarrow \quad G(t) \in \rL^q(dt)$$
	$$u(r) \in \rL^q(rdr) \quad \Longleftrightarrow \quad U(t) \in \rL^q(dt)$$
	$$f(r)/r^2 \in \rL^q(rdr) \quad \Longleftrightarrow \quad \psi(t) \in \rL^q(dt).$$
	Given $G(t)\in \rL^q(dt)$, our new goal is to solve for $U(t)\in \rL^q(dt)$ (or alternatively $\psi(t)\in W^{2,q}(dt)$).
	The ordinary differential equations in terms of the new variable and functions are
	$$U(t) = \psi''(t) + \bigg(4-\frac{4}{q}\bigg)\psi'(t) +\bigg(4-(km)^2+\frac{4}{q^2}-\frac{8}{q}\bigg)\psi(t)$$
	and
	$$\frac{1}{\alpha}U'(t)+(1-\frac{2}{\alpha q}-\lambda_1-i\lambda_2)U(t)-mik\beta e^{-\alpha t}U(t) -mik\alpha(2-\alpha)\beta e^{-\alpha t}\psi(t)=G(t).$$
	We distinguish two cases. In the case when $k=0$, we can integrate the equation of first order and get:
	$$U(t) = c e^{t(2/q+\alpha(\lambda-1))} + \alpha e^{t(2/q+\alpha(\lambda-1))}\int_0^t e^{-s(2/q+\alpha(\lambda-1))}G(s) ds.$$
	The unique choice of constant that will make $U$ integrable is precisely
	$$c= -\alpha \int_0^\infty e^{-s(2/q+\alpha(\lambda-1))}G(s)ds,$$
	with which we have
	$$U(t) = -\alpha\int_\mathbb{R} \chi_{(-\infty,0)}(t-s) e^{(t-s)(2/q+\alpha(\lambda-1))}G(s)ds.$$
	By Young's convolution inequality, it follows that $U(t)\in \rL^q(dt)$ is a solution of the ordinary differential equation with 
	$$\|U\|_q \leq \frac{\alpha \|G\|_{q}}{2/q+\alpha(\lambda_1-1)}.$$
	This is the unique solution, since in the case when $k=0$, the homogeneous problem reduces to :
	$$\frac{1}{\alpha}U'(t)+(1-\frac{2}{\alpha q}-\lambda_1-i\lambda_2)U(t)=0,$$
	whose only solution is given by $U(t) = c_1 e^{t(2/q+\alpha(\lambda-1))}$, which is not in any $\rL^q$ space, unless identically zero.

	We may henceforth without loss of generality assume that $k\geq 1$. The case when $k\leq -1$ is "symmetric" to the case when $k\geq 1$. Now, integrating the differential equation of first order yields
	$$U(t) = c_1\exp\big(-mik\beta e^{-\alpha t} +t(2/q+\alpha(\lambda-1))\big) +$$ $$+\alpha\exp\big(-mik\beta e^{-\alpha t} +t(2/q+\alpha(\lambda-1))\big)\int_0^t \exp\big(mik\beta e^{-\alpha s} -s(2/q+\alpha(\lambda-1))\big)G(s)ds+$$
	$$+mik\beta \alpha^2(2-\alpha)\exp\big(-mik\beta e^{-\alpha t} +t(2/q+\alpha(\lambda-1))\big)\times$$ \begin{equation}\label{INTEQN1}\times \int_0^t \exp\big(mik\beta e^{-\alpha s} -s(2/q+\alpha\lambda)\big)\psi(s)ds.\end{equation}
	The constant $c_1$ depends on the choice of initial condition. Also, integrating the differential equation of second order gets us: $$\psi(t) = c_2 e^{-(km+2-2/q)t} -\frac{e^{-(km+2-2/q)t}}{2km}\int_0^t e^{(km+2-2/q)s}U(s)ds + $$ \begin{equation} \label{INTEQN2}+c_3 e^{(km-2+2/q)t}+\frac{e^{(km-2+2/q)t}}{2km}\int_0^t e^{-(km-2+2/q)s} U(s)ds.\end{equation}
	
	Finally to have $\psi(t)\in \rL^q(dt)$ we use Equation \eqref{INTEQN2}, with the unique choice of constants:
	$$c_2 =  -\frac{1}{2km}\int_{-\infty}^0 e^{(km+2-2/q)s}U(s)ds$$
	and
	$$c_3= -\frac{1}{2km}\int_{0}^\infty e^{-(km-2+2/q)s}U(s)ds.$$
	These bounded linear functionals of $U$ are uniquely chosen so that $\psi \in \rL^q(dt)$ if $U\in \rL^q(dt)$. In addition, for smooth functions $u,f$, Equation \eqref{INTEQN2} implies that $\Delta^{-1}u = f$ in the classical sense, while the unique choice of $c_2,c_3$ above guarantee the desired integrability of $f$. We may thus write $\psi$ in terms of the convolution operator:
	$$\psi(t) = -\frac{1}{2km}\int_{\mathbb{R}} \bigg( e^{-(km+2-2/q)(t-s)}\chi_{(0,\infty)}(t-s) + e^{(km-2+2/q)(t-s)}\chi_{(-\infty, 0)}(t-s)\bigg)U(s)ds$$
	We can denote the kernel above by 
	$$K_1(t,s):= \bigg( e^{-(km+2-2/q)(t-s)}\chi_{(0,\infty)}(t-s) + e^{(km-2+2/q)(t-s)}\chi_{(-\infty, 0)}(t-s)\bigg).$$
	We remark that given a function $g$, there is also at most one unique choice of the parameter $c_1$ such that the function $U(t)$ lies in $\rL^q(dt)$ and satisfies Equation \eqref{INTEQN1}. In fact, the bounded functional $c_1(U)$ will be given by
	$$c_1 = -\alpha \int_0^\infty \exp\big(mik\beta e^{-\alpha s} -s(2/q+\alpha(\lambda-1))\big)G(s) ds-$$
	$$ mik\beta \alpha^2(2-\alpha)\int_0^\infty \exp\big(mik\beta e^{-\alpha s} -s(2/q+\alpha\lambda)\big)\psi(s)ds$$
	and Equation \eqref{INTEQN1} becomes
	$$U(t) = -\alpha \int_{\mathbb{R}} \exp\big(-mik\beta e^{-\alpha t} +mik\beta e^{-\alpha s} +(t-s)(2/q+\alpha(\lambda-1))\big)\chi_{(-\infty,0)}(t-s)G(s)ds - $$
	$$mik\beta \alpha^2(2-\alpha)\times $$ $$\times \int_{\mathbb{R}}e^{-\alpha s}\exp\big(-mik\beta e^{-\alpha t} +mik \beta e^{-\alpha s} +(t-s)(2/q+\alpha(\lambda-1))\big)\chi_{(-\infty,0)}(t-s)\psi(s)ds.$$
	We denote the kernel:
	$$K_2(t,s) := \exp\big(-mik\beta e^{-\alpha t} +mik\beta e^{-\alpha s} +(t-s)(2/q+\alpha(\lambda-1))\big)\chi_{(-\infty,0)}(t-s).$$
	We thus have the following Fredholm integral equation for $U$:
	$$U(t) = -\alpha \int_\mathbb{R} K_2(t,s)G(s)ds + \frac{i \alpha^2 (2-\alpha)\beta}{2}\int_{\mathbb{R}}\int_{\mathbb{R}} K_2(t,s)e^{-\alpha s}K_1(s,r)U(r)dr ds.$$
	
	We proceed to simplify the integral equation by integrating the kernels in the $s$ variable first:
	$$\int_\mathbb{R} K_2(t,s)e^{-\alpha s}K_1(s,r) ds = $$
	$$= F_1(t,r)\int_\mathbb{R}\exp\big(mik\beta e^{-\alpha s}-\alpha s -(2/q+\alpha(\lambda-1)+mk+2-2/q)s\big) \chi_{(\max(t,r),\infty)}(s)ds+ $$
	$$+F_2(t,r) \int_\mathbb{R}\exp\big(mik\beta e^{-\alpha s}-\alpha s -(2/q+\alpha(\lambda-1)-(mk-2+2/q))s\big) \chi_{(t,r)}(s)ds.$$
	Where 
	$$F_1(t,r) = \exp\big(-mik\beta e^{-\alpha t} + (2/q+\alpha(\lambda-1))t + r(mk+2-2/q)\big)$$
	and
	$$F_2(t,r) = \exp\big(-mik\beta e^{-\alpha t} + (2/q+\alpha(\lambda-1))t - r(mk-2+2/q)\big).$$
	
	We remark that
	$$\frac{d}{ds}\bigg( \frac{i}{m\alpha k \beta } e^{mik\beta e^{-\alpha s}}\bigg) = e^{mik\beta e^{-\alpha s} -\alpha s}.$$
	Thus integrating by parts and simplifying we get (for $\mu\in \mathbb{C}$ with $Re(\mu)>0$):
	\begin{equation}\label{NEATIDEN1}\int_t^\infty e^{mik\beta e^{-\alpha s} -\alpha s} e^{-\mu s} ds = \frac{-i}{m\alpha k\beta} e^{mik\beta e^{-\alpha t}}e^{-\mu t}+\frac{\mu i}{m\alpha k\beta} \int_t^\infty e^{mik\beta e^{-\alpha s}} e^{-\mu s}ds. \end{equation}
	Likewise, integrating by parts yields (for some generic $\mu\in \mathbb{C}$):
	$$\int_t^r e^{mik\beta e^{-\alpha s} -\alpha s} e^{-\mu s} ds =$$ \begin{equation}\label{NEATIDEN2} =\frac{-i}{m\alpha k \beta}\bigg(e^{mik\beta  e^{-\alpha t}}e^{-\mu t} -e^{mik\beta  e^{-\alpha r}}e^{-\mu r}\bigg) + \frac{\mu i}{m\alpha k\beta} \int_t^r e^{mik\beta e^{-\alpha s}}e^{-\mu s}ds.\end{equation}
	Henceforth we denote
	\begin{align*}
    \mu_1 &:= 2/q+a(\lambda-1)+mk+2-2/q,\\
    \mu_2 &:= 2/q+\alpha(\lambda-1)-(mk-2+2/q).
        \end{align*}
	
	Applying Equation \eqref{NEATIDEN1} and Equation \eqref{NEATIDEN2} to the original integral under consideration (and using that $\lambda>1-\frac{2}{q\alpha}$) gets us
	$$\int_\mathbb{R} K_2(t,s)e^{-\alpha s}K_1(s,r) ds = $$
	$$F_1(t,r)\chi_{(0,\infty)}(t-r)\bigg(\frac{-i}{m\alpha k\beta} e^{mik \beta e^{-\alpha t}}e^{-\mu_1 t} + \frac{\mu_1 i}{m\alpha k \beta}\int_t^\infty e^{mik\beta e^{-\alpha s}} e^{-\mu_1 s} ds\bigg)+$$
	$$+F_1(t,r)\chi_{(-\infty, 0 )}(t-r)\bigg(\frac{-i}{m\alpha k\beta} e^{mik \beta e^{-\alpha r}}e^{-\mu_1 r} + \frac{\mu_1 i}{m\alpha k \beta}\int_r^\infty e^{mik\beta e^{-\alpha s}} e^{-\mu_1 s} ds\bigg)+$$
	$$+F_2(t,r)\chi_{(-\infty,0)}(t-r)\bigg( \frac{-i}{m\alpha k \beta}\bigg(e^{mik\beta  e^{-\alpha t}}e^{-\mu_2 t} -e^{mik\beta  e^{-\alpha r}}e^{-\mu_2 r}\bigg) +$$ $$+ \frac{\mu_2 i}{m\alpha k\beta} \int_t^r e^{mik\beta e^{-\alpha s}}e^{-\mu_2 s}ds\bigg).$$

	Simplifying further (there is a remarkable cancellation that occurs) we achieve:
	$$ \int_\mathbb{R} K_2(t,s)e^{-\alpha s}K_1(s,r) ds = $$
	$$=\frac{-i}{m\alpha k \beta} K_1(t,r)+\frac{\mu_1 i}{m\alpha k \beta} \int_{\mathbb{R}} K_2(t,s) e^{-(s-r)(mk+2-2/q)}\chi_{(0,\infty)}(s-r)ds+$$
	$$+ \frac{\mu_2 i}{m\alpha k \beta} \int_\mathbb{R} K_2(t,s)e^{(s-r)(mk-2+2/q)}\chi_{(-\infty,0)}(s-r)ds.$$
	We rearrange the terms and get:
	$$ \int_\mathbb{R} K_2(t,s)e^{-\alpha s}K_1(s,r) ds = $$
	$$=\frac{-i}{m\alpha k\beta} K_1(t,r) + \frac{\mu_1 i}{m\alpha k \beta} \int_{\mathbb{R}} K_2(t,s)K_1(s,r) ds+$$
	$$+\frac{-2i}{\beta\alpha}\int_{\mathbb{R}}K_2(t,s) e^{(s-r)(mk-2+2/q)}\chi_{(-\infty,0)}(s-r)ds.$$
	Henceforth we denote
	$$K_3(s,r) := e^{(s-r)(mk-2+2/q)}\chi_{(-\infty,0)}(s-r).$$

	We conclude that
	$$U(t) = -\alpha \int_\mathbb{R} K_2(t,s)G(s)ds +\frac{\alpha(2-\alpha)}{2km}\int_\mathbb{R} K_1(t,r) U(r) dr +$$
	$$+\frac{-\mu_1\alpha(2-\alpha)}{2km}\iint_{\mathbb{R}^2} K_2(t,s)K_1(s,r) U(r) ds dr+$$ $$+\alpha(2-\alpha)\iint_{\mathbb{R}^2} K_2(t,s) K_3(s,r) U(r) ds dr.$$
	
	\subsection{The Integral Transforms}
	We have a convolution operator:
	$$\Phi_1(U)(t) :=  \int_{\mathbb{R}} K_1(t,s) U(s)ds.$$
	Now, by Young's convolution inequality, we have
	$$\|\Phi_1(U)(t)\|_q \leq \|U\|_q \bigg(\int_{0}^\infty e^{-(mk+2-2/q)y} dy+\int_{-\infty}^0 e^{(mk-2+2/q)y}dy\bigg)\leq$$ $$\leq \bigg(\frac{1}{mk-2+2/q}+\frac{1}{mk+2-2/q}\bigg) \|U\|_q.$$
	In a similar manner, we denote
	$$\Phi_2(U):=\int_{\mathbb{R}}K_2(t,s) U(s)ds$$
	and observe that
	$$\abs{\Phi_2(U)} \leq \int_{\mathbb{R}}\hat{K_2}(t,s)\abs{U(s)}ds$$
	where
	$$\hat{K_2}(t,s):= e^{(t-s)(2/q+\alpha(\lambda-1))}\chi_{(-\infty,0)}(t-s)$$
	is the kernel of a convolution operator. Thus, arguing as before, we apply the Young convolution inequality and get
	$$\|\Phi_2(U)(t)\|_q \leq \|U\|_q \bigg(\int_{0}^\infty e^{-(2/q+\alpha(\lambda_1-1))y} dy\bigg)\leq \frac{1}{2/q+\alpha(\lambda_1-1)} \|U\|_q.$$
	Lastly, we have the convolution operator:
	$$\Phi_3(U)(t) :=  \int_{\mathbb{R}} K_3(t,s) U(s)ds.$$
	whose operator norm we once again estimate using Young's convolution inequality:
	$$\|\Phi_3(U)(t)\|_q \leq \|U\|_q \bigg(\int_{-\infty}^0 e^{y(mk-2+2/q)} dy\bigg)\leq \frac{1}{(mk-2+2/q)} \|U\|_q.$$

	\subsection{The Fixed Point Argument}
	We recall what we have heretofore shown. We found three bounded integral transforms from $\rL^q\to \rL^q$, given by $\Phi_1$, $\Phi_2$ and $\Phi_3$. Showing that the operator we are working on is surjective (and injective) then becomes equivalent to showing (for every $G\in \rL^q$) the existence of a (unique) fixed point of the map 
	$$U \to -\alpha \Phi_2(G) +\frac{\alpha(2-\alpha)}{2km} \Phi_1(U) - \frac{\mu_1 \alpha(2-\alpha)}{2km}\Phi_2(\Phi_1(U))+\alpha(2-\alpha)\Phi_2(\Phi_3(U)).$$
	This, by using the Banach fixed point theorem (or a Neumann series), is equivalent to showing that the map 
	$$\hat{\Phi}(U):=\frac{\alpha(2-\alpha)}{2km} \Phi_1(U) - \frac{\mu_1 \alpha(2-\alpha)}{2km}\Phi_2(\Phi_1(U))+\alpha(2-\alpha)\Phi_2(\Phi_3(U))$$ 
	is contractive, i.e. $\|\hat{\Phi}\|<1$. However, by our work in the previous subsection, we can estimate the operator norm of $\hat{\Phi}$ by:
	$$\|\hat{\Phi}\|\leq\frac{\alpha(2-\alpha)}{2km}\cdot \bigg(\frac{1}{mk-2+2/q}+\frac{1}{mk+2-2/q}\bigg) +$$
	$$+ \frac{(2/q+a(\lambda_1-1)+mk+2-2/q) \alpha(2-\alpha)}{2km(2/q+\alpha(\lambda_1-1))}\cdot\bigg(\frac{1}{mk-2+2/q}+\frac{1}{mk+2-2/q}\bigg)+$$
	$$+ \frac{\alpha(2-\alpha)}{2/q+\alpha(\lambda_1-1)}\bigg( \frac{1}{mk-2+2/q}\bigg).$$
	Recall that the complex number $\lambda$ was chosen so that $\textrm{Re}(\lambda) = \lambda_1> 1-\frac{2}{q\alpha}$. Let $\epsilon>0$ be such that 
	$$\lambda_1 = 1-\frac{2}{q\alpha} +\epsilon = 1-\frac{2-\epsilon q\alpha}{q\alpha}.$$
	Using Mathematica, we can verify that if
	$$mk \geq 2+\frac{8}{\epsilon \alpha q},$$
	then $\|\hat{\Phi}\|<1$ for all choices of $1< q\leq 2/\alpha$ and $0<\alpha<1$. For the code that verifies this inequality, see the last subsection of the appendix. Thus the unique solubility of the ordinary differential equation, hence the surjectivity of the map, is proven for the spaces $\rL^2_m$ or $U_{k,q}$ whenever $m\geq 2+\frac{8}{\epsilon\alpha q }$ or $k\geq 2+\frac{8}{\epsilon \alpha q}$ and $1< q \leq 2/\alpha$,  $0<\alpha<1$. Note that, since the fixed point is unique, we have provided a proof of the injectivity of the map $L_{ss}-\lambda I$ for $m, \lambda, k$ satisfying the conditions above.
	
	We now conclude. Observe that $\epsilon = \lambda_1 -(1-\frac{2}{q\alpha})$, then we have proven that $\lambda\in \rho_\alpha(L_{ss})$ if 
	$$\textrm{Re}(\lambda) \geq 1-\frac{2}{\alpha q} +\frac{8}{\alpha q(m-2)},$$
	which directly implies the upper estimate on the spectral bound.

    \subsection{Lower Bound for the Spectral Bound}
	\label{ssec:spectral bound}
	
	In this subsection we show that  the line $\textrm{Re}(\lambda)=1-1/\alpha$ is contained in the spectrum of $L_{ss}$ when defined on $\rL^2_m$. More precisely, we have the following 
	
	\begin{prop}
		$\lambda = 1-\frac{1}{\alpha}-i \gamma\in \sigma_{\alpha}(L_{ss},\rL^2_m)$ for all $\alpha\in (0,1)$, $m\geq 2$, and $\gamma\in \mathbb{R}$. 
	\end{prop}

	\begin{proof}  
		From \autoref{ssec:prelims} we know that 
		it suffices to show that $1-\frac{1}{\alpha}+i\gamma\in\sigma_{\alpha}(L_{ss},U_{km})$ for some $k\in \mathbb{Z}$. We do this for $k=0$, in which case the resolvent equation for $U(t)$, as above, simplifies to:
		\begin{equation}\label{DERIVBAD}
			U'(t)+i\alpha\gamma U(t)=\alpha G(t).
		\end{equation}
		Integrating, Equation \eqref{DERIVBAD} becomes the integral equation
		$$U(t) = ce^{-i\alpha\gamma t} + \alpha e^{-i\alpha\gamma t}\int_0^t e^{i\alpha \gamma s}G(s) ds,$$
		where $c$ is a constant that must be chosen so that $U\in \rL^2(\mathbb{R})$. Consider the function
		$$G(t) := \frac{e^{-i\alpha\gamma s}}{\abs{t}+1}\in \rL^2(\mathbb{R}).$$
		It is not difficult to see that 
		$$ce^{-i\alpha\gamma t} + \alpha e^{-i\alpha\gamma t}\int_0^t e^{i\alpha \gamma s}G(s) ds = e^{-i\alpha\gamma t}\big(c + \alpha\cdot\textrm{sgn}(t)\log(1+\abs{t})\big),$$
		is not in $\rL^2(\mathbb{R})$ for any choice of the constant $c$. This implies that $1-\frac{1}{\alpha}-i\gamma \in \sigma_\alpha(L_{ss}, U_0)$ for all $\gamma\in \mathbb{R}$.
	\end{proof}
	
	As a direct consequence, we obtain for the spectral bound of $L_{ss}$ the following result.
	
	\begin{cor}\label{cor:spectral bound}
		The spectral bound of $L_{ss}$ is bounded from below by
		\begin{equation*}
			1-\frac{1}{\alpha} \leq s(L_{ss}) .
		\end{equation*}
	\end{cor}

	\section{Explicit Solutions of the Self-Similar Eigenvalue Equation}
	\label{sec:explicit solution}
	
		In this section we solve the eigenvalue equation in terms of the stream function.
	For simplicity of notation, we assume $\beta=1$ and $m=2$, but the proof is entirely the same in the general case. We recall Equation \eqref{HOMOGODE} :
	
	$$
	\frac{r}{\alpha}\big(f_k''(r)+\frac{f_k'(r)}{r}-\frac{4k^2f_k(r)}{r^2}\big)'+(1-(\lambda_1+i\lambda_2))\big(f_k''(r)+\frac{f_k'(r)}{r}-\frac{4k^2f_k(r)}{r^2}\big)$$
	$$
	-2ikr^{-\alpha}\big(f_k''(r)+\frac{f_k'(r)}{r}-\frac{4k^2f_k(r)}{r^2}\big)-2ikr^{-2-\alpha}\alpha(2-\alpha)f_k(r)=0.
	$$
	
	We consider the above homogeneous ordinary differential equation and perform a change of variables $z= -2ik r^{-\alpha}$, since this will put the ordinary differential equation in a form that is already well studied. More precisely, we define a new function $w_k(z)$ that is given by
	$$w_k(z) = f_k((-2ik)^{1/\alpha} z^{-1/\alpha})\cdot (-2ik)^{-2k/\alpha}z^{2k/\alpha}.$$
	The function is defined so that if $f_k(r) = w(-2ik r^{-\alpha}) r^{2k}$ satisfies the homogeneous ordinary differential equation, we get that $w$ satisfies the following ordinary differential equation:
	
	$$0=2 k r^{-3 \alpha+2 k-2}\bigg(i (\alpha-4 k) r^\alpha (r^\alpha (\alpha \lambda-2 k+2)+2 i \alpha k) w'(-2 i
	k r^{-\alpha})+$$ $$+2 \alpha k ((r^\alpha (\alpha (\lambda+2)-6 k+2)+2 i \alpha k) w''(-2 i k
	r^{-\alpha})-2 i \alpha k w^{(3)}(-2 i k r^{-\alpha}))+$$ $$+i (\alpha-2) \alpha r^{2 \alpha}
	w(-2 i k r^{-\alpha})\bigg)$$
	
	Letting $r=(-2ik)^{1/\alpha} z^{-1/\alpha}$ and simplifying we finally have:
	$$0=  z^2 w^{(3)}(z)+z  (1-z+\lambda+1+\frac{2-6k}{\alpha})w''(z)+$$ \begin{equation}\label{HOMO2}+
		((\lambda+\frac{2-2k}{\alpha})\alpha^{-1}(\alpha-4 k)-z(-1+\alpha^{-1}(\alpha-4 k)+1))w'(z)+\alpha^{-1}(\alpha-2)  w(z).\end{equation}
	
	We put the equation in this form because Equation 07.25.13.0004.01 from \cite{W} states that the general solution of the ordinary differential equation:
	$$z^2w^{(3)}(z)+z(1-z+b_1+b_2)w''(z)+(b_1b_2-z(a_1+a_2+1))w'(z)-a_1a_2w(z)=0$$ is given exactly by 
	$$w(z) = c_1\cdot  \, _2\tilde{F}_2(a_1, a_2 ; b_1, b_2 ; z) +$$
	$$+c_2\bigg( G^{2,2}_{2,3}\big(z \big| 1-a_1, 1-a_2 ; 0, 1-b_1, 1-b_2\big) +G^{2,2}_{2,3}\big(z \big| 1-a_1, 1-a_2 ; 0, 1-b_2, 1-b_1\big) \bigg) +$$ $$+ c_3 G^{3,2}_{2,3}\big(-z \big|  1-a_1, 1-a_2 ; 0, 1-b_1, 1-b_2\big).$$
	
	Here $\, _2\tilde{F}_2$ denotes the regularized hypergeometric function and $G^{p,q}_{m,n}$ denotes the Meijer-G function. The definition of both of these classes of special functions can be found in \cite{LUKE} or \cite{W}.
	
	We choose 
	$$a_1 = \frac{-2k-\sqrt{\alpha^2-2\alpha+4k^2}}{\alpha}, \quad a_2 = \frac{-2k+\sqrt{\alpha^2-2\alpha+4k^2}}{\alpha}$$
	$$b_1 = \frac{\alpha-4k}{\alpha}, \quad b_2 = \frac{2-2k+\alpha\lambda}{\alpha}$$
	and denote
	$$\mathfrak{q}_{\alpha,k}:=\frac{\sqrt{\alpha^2-2\alpha+4k^2}}{\alpha}.$$
	Then the solution of Equation \eqref{HOMO2} is exactly (for any choice of constants $c_1,c_2,c_3$):
	$$w(z) = c_1\cdot  \, _2\tilde{F}_2(-\frac{2k}{\alpha}-\mathfrak{q}_{\alpha,k}, -\frac{2k}{\alpha}+\mathfrak{q}_{\alpha,k} ; 1-\frac{4k}{\alpha}, \lambda+\frac{2-2k}{\alpha} ; z) +$$
	$$+c_2\bigg( G^{2,2}_{2,3}\big(z \big| 1+\frac{2k}{\alpha}+\mathfrak{q}_{\alpha,k}, 1+\frac{2k}{\alpha}-\mathfrak{q}_{\alpha,k} ; 0, \frac{4k}{\alpha}, 1-\lambda+\frac{2k-2}{\alpha}\big) +$$ $$+G^{2,2}_{2,3}\big(z \big| 1+\frac{2k}{\alpha}+\mathfrak{q}_{\alpha,k}, 1+\frac{2k}{\alpha}-\mathfrak{q}_{\alpha,k} ; 0,1-\lambda+\frac{2k-2}{\alpha}, \frac{4k}{\alpha}\big)  \bigg) +$$ $$+ c_3 G^{3,2}_{2,3}\big(-z \big|  1+\frac{2k}{\alpha}+\mathfrak{q}_{\alpha,k}, 1+\frac{2k}{\alpha}-\mathfrak{q}_{\alpha,k} ; 0, \frac{4k}{\alpha}, 1-\lambda+\frac{2k-2}{\alpha}\big).$$
	Using the transformation back to $f_k(r)$, we see that we have found the solutions of the homogeneous differential equation \eqref{HOMOGODE}. The asymptotic properties for the special functions $\, _2\tilde{F}_2$ and $G^{p,q}_{m,n}$ described in \cite[Chapter 5]{LUKE}, show that no solution of the homogeneous equation is integrable, which implies that the point spectrum of $L_{ss}$ is empty when considered on $L^2_m$, indepndently of $m$.

	\section{Mathematica Code}
	\label{sec:mathematica}
	
		The following Mathematica code verifies that the previously discussed contraction property is satisfied for $m\geq 2+\frac{8}{\epsilon q\alpha}$. We should interpret $e$ in the code below as representing the number $\epsilon q\alpha$ from our proof above.
	
	\begin{lstlisting}
		Reduce[ForAll[{a, l, k, q}, 
		0 < a < 1 && k >= 1 && l >= 1 - (2 - e)/(q*a) && 
		1 <= q <= 
		2/a, (a*(2 - a)/(2*k*m))*(1/(k*m + 2 - 2/q) + 
		1/(k*m - 2 + 2/q)) + ((2/q + a*(l - 1) + 
		m*k + 2 - 2/q)*
		a*(2 - a)/(2*k*m*(2/q + a*(l - 1))))*
		(1/(k*m + 2 - 2/q)
		+ 
		1/(k*m - 2 + 2/q)) + 
		a*(2 - a)*(1/(2/q + a (l - 1)))*
		(1/(m*k - 2 + 2/q)) < 1] && 
		e > 0 && m >= (8 + 2 e)/e && Element[m, Integers]]
	\end{lstlisting}

\section{\texorpdfstring{Linear Stability of the Power-Law Vortex in Physical Coordinates in $\rL^2(\R^2)$}{Linear Stability of the Power-Law Vortex in Physical Coordinates in rL2(R2)}}
\label{sec:physicalunweighted}

In this section we consider the linearization of the Euler equation around the power-law vortex in physical coordinates in $\rL^2(\R^2)$. The goal is to prove \autoref{prop:physicalMAINMAINunweighted}:

	\begin{prop}\label{prop:1}
		Let $\beta\in \mathbb{R}\setminus\{0\}$. Let $m\geq 1$ and $\alpha \in (0,1)$. % or $m=1$ and $\alpha \in (0,2-\sqrt{2})$. 
		Let $\overline{\Omega}(\xi) = \beta(2-\alpha) \abs{\xi}^{-\alpha}$ be the radial power-law vortex that solves the unforced Euler equations in self-similar coordinates and let $L_{\phi}$ be the linearization of the Euler equations around $\overline{\Omega}$. If $L_{\phi}$ generates a $C_0$-semigroup on $\rL^2_{m}(\R^2)$ then
		$$\sigma(L_\phi) \subset \{z \in \mathbb{C} : \mathrm{Re}(z)=0\}.$$
	\end{prop}
	
We now introduce the following invertible operator
$$S_s :\rL^2(\R^2) \to \rL^2(\R^2) \colon  \Omega(x) = \Omega(s\cdot x).$$
The inverse is evidently given by $S_s^{-1} = S_{s^{-1}}$. Note that 
\begin{equation}
	\label{DELCOMMUTE}
	\partial_i S_s U = s \partial_i \nabla U,
\end{equation}
where $\partial_i$ denotes the partial derivative along $e_i$. This implies that a similar formula for the gradient, the divergence and directionals derivatives, etc. and $U$ may denote a function, vector field, tensor, etc.

%Now let us suppose that we have a homogeneous self-similar solution $\Omega$. Evidently, $\Omega$ is a stationary solution of the Euler equations; therefore, we may linearize the Euler equations around $\Omega$ and obtain the operator 
%\begin{equation}
%	\label{LINEULER}
%	\begin{gathered}
%		L_\phi\Omega = -(\overline{U}\cdot \nabla)\Omega-(U\cdot \nabla)\overline{\Omega} +(\overline{\Omega}\cdot\nabla)U+(\Omega\cdot\nabla)\overline{U}\\
%		\mathrm{K}_{BS}\ast \Omega(\cdot) = U(\cdot)\\
%		\mathrm{K}_{BS}\ast \overline{\Omega}(\cdot) = \overline{U}(\cdot)
%	\end{gathered}
%\end{equation}
On the other hand we obtain the following commutative rule of $S$ with the Biot-Savart operator.
\begin{lem}
	\label{BIOTCOMMUTE}
	The Biot Savart kernel $K_{BS}$ has the property that
	$$K_{BS}\ast(S_s\Omega)= s^{-1}S_s(K_{BS}\ast \Omega).$$  
\end{lem}
\begin{proof}
	In $\mathbb{R}^n$, the Biot-Savart kernel, or the kernel of the operator $\nabla^\perp \Delta^{-1}$, is given, up to constants, by
	$$K_{BS}(y) = \frac{y^\perp}{|y|^n}.$$ Therefore, up to constants, we have
	$$K_{BS}\ast(S_s \Omega)(x) = \int_{\mathbb{R}^n} \frac{(x-y)^\perp}{|x-y|^n}\Omega(s y) dy.$$
	Performing the change of variables $s y\to y$, we get
	$$K_{BS}\ast(S_s \Omega)(x) = \int_{\mathbb{R}^n} \frac{s^{-1}(s x- y)^\perp}{|s x - y|^n} \Omega(y)dy= s^{-1}S_s(K_{BS}\ast \Omega). $$
    \qedhere
\end{proof}

Now, we are sufficient prepared to prove the \autoref{prop:1}. 

\begin{proof}[Proof of {\autoref{prop:1}}]
	Since $S_s$ is a bounded operator with bounded inverse, we obtain, for any $s>0$
	$$\sigma(L_\phi) = \sigma(S_s^{-1}L_\phi S_s).$$
	On the other hand, we may use Equation \eqref{DELCOMMUTE}, \autoref{BIOTCOMMUTE}
	and the fact that $\bar{\Omega}$ is $\alpha$-homogeneous 
	to conclude 
	$$S_s^{-1}L_\phi S_s = s^\alpha L_\phi.$$
	Combining both we obtain
	\begin{equation}\label{ESSEQ}
		\sigma(L_\phi) 
		= \sigma(S_s^{-1}L_\phi S_s) = \sigma( s^\alpha L_\phi) = s^\alpha \sigma(L_\phi).
	\end{equation}
	for any $s>0$.
	
	\smallskip 
	
	By assumption $L_\phi$ generates a $C_0$-semigroup. Since the Euler equation is time reversible we conclude that $-L_\phi$ does so was well. By \autoref{thm:generation groups} we conclude that $L_\phi$ generates a $C_0$-group. This implies that $\lambda \in \rho(L_{\phi})$ for $|\Re \lambda \, |$ sufficient large. 
	We conclude from \eqref{ESSEQ} that $\Re \lambda = 0$ and hence the claim.
\end{proof}
	
%	\section{Exponential Linear Stability in Self-Similar Coordinates}
%	
%	
%
%	An immediate consequence of \cite[Theorem 1.3]{C} is
%	\begin{lem}
%		\label{SURJECTIVEOP}
%		For $m\geq 3$ and for all $k\in \mathbb{Z}$, there exists $\lambda_0>0$ such that $\lambda_0\in \rho(C_0)$.
%	\end{lem}
%	\begin{proof}
%		Choose $\lambda_0 > \tfrac{4}{\alpha}$.
%	\end{proof}

\end{document}